\documentclass[11pt,leqno]{article}
\usepackage[doc]{optional}
\usepackage{color}
\usepackage{graphicx}
\definecolor{labelkey}{rgb}{0,0.08,0.45}
\definecolor{refkey}{rgb}{0,0.6,0.0}
\definecolor{Brown}{rgb}{0.45,0.0,0.05}
\usepackage{mathpazo}
\usepackage{amsmath}
\usepackage{amssymb}
\usepackage{theorem}
\newcommand{\fDK}{$\fD{}\,$-Klee}
\newcommand{\nnn}{\ensuremath{{n\in{\mathbb N}}}}
\newcommand{\thalb}{\ensuremath{\tfrac{1}{2}}}
\newcommand{\menge}[2]{\big\{{#1}~\big |~{#2}\big\}}

\newcommand{\To}{\ensuremath{\rightrightarrows}}

\newcommand{\fenv}[1]%
{\ensuremath{\,\overrightarrow{\operatorname{env}}_{#1}}}
\newcommand{\benv}[1]%
{\ensuremath{\,\overleftarrow{\operatorname{env}}_{#1}}}

\newcommand{\scal}[2]{\left\langle{#1},{#2}  \right\rangle}

\newcommand{\RR}{\ensuremath{\mathbb R}}
\newcommand{\RP}{\ensuremath{\left[0,+\infty\right[}}
\newcommand{\RPX}{\ensuremath{\left[0,+\infty\right]}}

\newcommand{\RX}{\ensuremath{\,\left]-\infty,+\infty\right]}}

\newcommand{\IDD}{\ensuremath{U}}

\newcommand{\dom}{\ensuremath{\operatorname{dom}}}
\newcommand{\argmin}{\ensuremath{\operatorname{argmin}}}
\newcommand{\argmax}{\ensuremath{\operatorname*{argmax}}}

\newcommand{\gr}{\ensuremath{\operatorname{gr}}}

\newcommand{\intdom}{\ensuremath{\operatorname{int}\operatorname{dom}}\,}

\newcommand{\bd}{\ensuremath{\operatorname{bdry}}}
\newcommand{\ran}{\ensuremath{\operatorname{ran}}}
\newcommand{\conv}{\ensuremath{\operatorname{conv}}}
\newcommand{\cconv}{\ensuremath{\overline{\operatorname{conv}}\,}}

\newcommand{\ffproj}[1]{\overrightarrow{Q\thinspace}_%
{\negthinspace\negthinspace #1}}

\newcommand{\fD}[1]{\overrightarrow{D\thinspace}_%
{\negthinspace\negthinspace #1}}
\newcommand{\ffD}[1]{\overrightarrow{F\thinspace}_%
{\negthinspace\negthinspace #1}}
\newcommand{\ffDbz}{\protect\overrightarrow{F\thinspace}_%
{\negthinspace\negthinspace C}(\bz)}

\newcommand{\minf}{\ensuremath{-\infty}}
\newcommand{\pinf}{\ensuremath{+\infty}}
\newcommand{\bx}{\ensuremath{\mathbf{x}}}
\newcommand{\bc}{\ensuremath{\mathbf{c}}}
\newcommand{\by}{\ensuremath{\mathbf{y}}}
\newcommand{\bz}{\ensuremath{\mathbf{z}}}
\newcommand{\bg}{\ensuremath{\mathbf{g}}}
\newcommand{\bh}{\ensuremath{\mathbf{h}}}

{\begin{list}{}{%
\settowidth{\labelwidth}{\textrm{#1~}}%
\setlength{\leftmargin}{\labelwidth+\labelsep}}}
{\end{list}}
\newtheorem{theorem}{Theorem}[section]
\newtheorem{lemma}[theorem]{Lemma}
\newtheorem{corollary}[theorem]{Corollary}
\newtheorem{proposition}[theorem]{Proposition}
\newtheorem{definition}[theorem]{Definition}
\theoremstyle{plain}{\theorembodyfont{\rmfamily}
}
\theoremstyle{plain}{\theorembodyfont{\rmfamily}
}
\theoremstyle{plain}{\theorembodyfont{\rmfamily}
}
\theoremstyle{plain}{\theorembodyfont{\rmfamily}
\newtheorem{example}[theorem]{Example}}
\newtheorem{fact}[theorem]{Fact}
\theoremstyle{plain}{\theorembodyfont{\rmfamily}
\newtheorem{remark}[theorem]{Remark}}

\newcommand{\boxedeqn}[1]{%
    \[\fbox{%
        \addtolength{\linewidth}{-2\fboxsep}%
        \addtolength{\linewidth}{-2\fboxrule}%
        \begin{minipage}{\linewidth}%
        \begin{equation}#1\\[+4mm]\end{equation}%
        \end{minipage}%
      }\]%
  }

\begin{document}

\title{\textrm{Klee sets and Chebyshev centers\\ for the right Bregman
distance}}

\author{
Heinz H.\ Bauschke\thanks{Mathematics, Irving K.\ Barber School,
UBC Okanagan, Kelowna, British Columbia V1V 1V7, Canada. E-mail:
\texttt{heinz.bauschke@ubc.ca}.},
~Mason S.\ Macklem\thanks{Mathematics, Irving K.\ Barber School, UBC Okanagan,
Kelowna, British Columbia V1V 1V7, Canada. E-mail:
\texttt{mason.macklem@ubc.ca}.},
~Jason B.\  Sewell\thanks{Mathematics, Irving K.\ Barber School, UBC Okanagan,
Kelowna, British Columbia V1V 1V7, Canada. E-mail:
\texttt{js1868pp@interchange.ubc.ca}.},
~and Xianfu\
Wang\thanks{Mathematics, Irving K.\ Barber School, UBC Okanagan,
Kelowna, British Columbia V1V 1V7, Canada.
E-mail:  \texttt{shawn.wang@ubc.ca}.}}

\date{August 11, 2009}
\maketitle

\vskip 8mm

\begin{abstract} \noindent
We systematically investigate the farthest distance function,
farthest points, Klee sets, and Chebyshev centers,
with respect to Bregman distances induced by
Legendre functions.
These objects are of considerable interest in Information Geometry
and Machine Learning; when
the Legendre function is specialized
to the energy, one obtains classical notions
from Approximation Theory and Convex Analysis.

The contribution of this paper is twofold.
First, we provide an affirmative answer to a
recently-posed question on
whether or not every Klee set with respect to the
right Bregman distance is a singleton.
Second, we prove uniqueness of the Chebyshev center
and we present a characterization that relates
to previous works by Garkavi, by Klee, and by Nielsen and Nock.
\end{abstract}

{\small


}

\section{Introduction}

Throughout this paper,
\boxedeqn{
\text{$\RR^J$ is
the standard Euclidean space with inner
product $\scal{\cdot}{\cdot}$ and induced norm $\|\cdot\|$.}
}

Suppose that $S$ is a nonempty subset of $\RR^J$ such that
for every point in $\RR^J$, there exists a unique farthest point in $S$,
where ``farthest'' is understood in the standard Euclidean distance sense.
Then $S$ is said to be a \emph{Klee set}, and it is known that
$S$ must be a \emph{singleton}; see, e.g.,
\cite{Edgar1,urruty2,urruty3,Klee61,MSV} for further information.
(The situation in Hilbert space remains unsettled to this day.)

In \cite{bwyy2}, Klee sets were revisited from a new perspective by
using measures of fairly different from distances induced by norms.
To describe and follow up on
this viewpoint, we assume throughout that
\boxedeqn{
\label{eq:preD}
\text{$f \colon \RR^J\to\RX$
is a convex function of Legendre type.}
}
Recall that for a convex function
$g\colon\RR^J\to\RX$, the \emph{(essential) domain} is
$\dom g = \menge{x\in\RR^J}{g(x)\in \RR}$ and
$x^*\in\RR^J$ is a \emph{subgradient} of $g$ at a point $x\in\dom g$,
written $x^*\in\partial g(x)$, if $(\forall h\in\RR^J)$
$g(x) + \scal{h}{x^*} \leq g(x+h)$; this induces the corresponding
set-valued \emph{subdifferential operator} $\partial g\colon\RR^J\To\RR^J$.
(For basic terminology and results from
Convex Analysis not stated explicitly in this paper,
we refer the reader to \cite{lewis,Rock70,Simons,Zali02}.)
Then $g$ is said to be
\emph{essentially smooth} if $g$ is differentiable on $\intdom g$
 (the interior of its domain),
and $\|\nabla g(x)\|\to\pinf$ whenever $x$ approaches a point in the
boundary $\bd\dom g$;
$g$ is \emph{essentially strictly convex} if $g$ is strictly convex on
every convex subset of $\dom\partial g = \menge{x\in\RR^J}{\partial
g(x)\neq\varnothing}$; and
$g$ is a \emph{convex function of Legendre type} --- often simply
called a \emph{Legendre function} --- if $g$ is both essentially smooth and
essentially strictly convex.
See \cite{commun01,BorVan01,BorVanBook,Rock70} for further information on
Legendre functions.
It will be convenient to set
\boxedeqn{ U := \intdom f.
}

Many examples of Legendre functions exist; however, in this paper,
we focus mainly on the following.

\begin{example}[Legendre functions] \label{ex:Legendre}
The following are Legendre functions, each
evaluated at a point $x\in\RR^J$.
\begin{enumerate}
\item
\label{ex:Legendre:i}
\emph{(halved) energy:}
$f(x) = \thalb\|x\|^2 = \thalb\sum_{j}x_j^2$.
\item
\label{ex:elliptic:iv}
$f(x) = \thalb\scal{x}{Ax}$, where $A\in\RR^{J\times J}$ is symmetric and
positive definite.
\item
\label{ex:Legendre:ii}
\emph{negative entropy:}
$f(x)= \displaystyle \begin{cases}
\sum_{j} \big(x_j\ln(x_j)-x_j\big), &\text{if $x\in\RR^J_+$;}\\
\pinf, &\text{otherwise.}
\end{cases}$
\item
\label{ex:Legendre:iii}
\emph{negative logarithm:}
$f(x)= \displaystyle \begin{cases}
-\sum_{j} \ln(x_j), &\text{if $x\in\RR^J_{++}$;}\\
\pinf, &\text{otherwise.}
\end{cases}$
\end{enumerate}
Note that $U=\RR^J$ in \ref{ex:Legendre:i} and \ref{ex:elliptic:iv},
whereas $U=\RR^J_{++}$ in \ref{ex:Legendre:ii}
and \ref{ex:Legendre:iii}.
\end{example}

Legendre functions are of considerable interest to us because they
give rise to a very nice measure of discrepancy between points,
nowadays termed the ``Bregman distance''; see, e.g.,
\cite{Bregman,ButIus,CenZen}.

\begin{definition}[Bregman distance]
The \emph{Bregman distance} with respect to $f$, written $D_f$ or
simply $D$, is the function
\boxedeqn{\label{eq:D}
D\colon \RR^J \times \RR^J \to \RPX \colon (x,y) \mapsto
\begin{cases}
f(x)-f(y)-\scal{\nabla f(y)}{x-y}, &\text{if}\;\;y\in\IDD;\\
\pinf, & \text{otherwise}.
\end{cases}
}
\end{definition}

Although well established,
the term ``Bregman distance'' is a misnomer because
a Bregman distance is in general neither symmetric nor does it satisfy
the triangle inequality.
However, the Bregman distance is able to distinguish
different points in the
sense (see \cite[Theorem~3.7.(iv)]{Baus97}) that
\begin{equation}
\label{e:misnomer}
\big(\forall x\in \RR^J\big)
\big(\forall y\in U\big)
\quad
D(x,y)=0 \;\;
\Leftrightarrow\;\;
x=y.
\end{equation}

\begin{example}[Bregman distances]
\label{ex:distances}
The Bregman distances corresponding to the Legendre functions of
Example~\ref{ex:Legendre} between two points $x$ and $y$ in $\RR^J$
are as follows.
\
\begin{enumerate}
\item
\label{ex:distances:E}
\emph{(halved) Euclidean distance squared:}
$D(x,y)= \thalb\|x-y\|^2$.
\item
\emph{(halved) Mahalanobis distance squared:}
$D(x,y) = \thalb\scal{x-y}{A(x-y)}$.
\item
\label{ex:distances:KL}
\emph{Kullback-Leibler divergence:}
$D(x,y)= \displaystyle \begin{cases}
\sum_{j} \big(x_j\ln(x_j/y_j)-x_j+y_j\big), &\text{if $x\in\RR^J_+$ and
$y\in\RR^J_{++}$;}\\
\pinf, &\text{otherwise.}
\end{cases}$
\item
\label{ex:distances:IS}
\emph{Itakura-Saito distance:}
$D(x,y)= \displaystyle \begin{cases}
\sum_{j} \big(\ln(y_j/x_j)+x_j/y_j - 1\big), &\text{if $x\in\RR^J_{++}$ and
$y\in\RR^J_{++}$;}\\
\pinf, &\text{otherwise.}
\end{cases}$
\opt{hhb}{{\color{blue}\marginpar{hhb}
\\$\blacktriangleright$
Observe that this one is not strictly convex (even not convex) in $y$!
$\blacktriangleleft$
}}
\end{enumerate}
\end{example}

From now on, we assume that $C$ is a subset of $\RR^J$ such that
\boxedeqn{
\label{e:CsubU}
\varnothing\neq C\subseteq U.
}

\begin{definition}[right Bregman farthest-distance function and
farthest-point map]
The \emph{right Bregman farthest-distance function} is
\boxedeqn{\label{e:montag:a}
\ffD{C} \colon \RR^J \to \RPX \colon  x\mapsto \sup_{c\in C}D(x,c),
}
and the corresponding \emph{right Bregman farthest-point map} is
\boxedeqn{ \label{e:montag:b}
\ffproj{C} \colon\RR^J \To\RR^J \colon
x\mapsto
\begin{cases}
\argmax_{c\in C} D(x,c) = \menge{c\in C}{D(x,c) = \ffD{C}(x)}, &\text{if $x\in\dom f$;}\\
\varnothing, &\text{otherwise.}
\end{cases}
}
\end{definition}

Observe that
\begin{equation}
\label{e:090217:a}
\ffD{C} \;\;\text{is convex and lower semicontinuous,}
\end{equation}
\opt{hhb}{{\color{blue}\marginpar{\color{blue}hhb}
$\blacktriangleright$ because it is the supremum of such functions $x\mapsto
f(x)-f(c)-\scal{x-c}{\nabla f(c)}$
$\blacktriangleleft$
}}
that
\begin{equation}
\label{e:090216:f}
\dom \ffproj{C} \subseteq \dom \ffD{C} \subseteq \dom f,
\end{equation}
and that
\begin{equation}
\label{e:090322n}
\text{if $C$ is compact, then }\dom\ffproj{C} = \dom \ffD{C} =
\dom f.
\end{equation}
We are now ready to continue the discussion on Klee sets started earlier
by introducing a notion central to this paper.

\begin{definition}[\fDK\ set]
The set $C$ is said to be \fDK, if
for every $x\in U$, $\ffproj{C}x$ is a singleton.
\end{definition}

The asymmetry of $D$ gives also rise to
the \emph{left} Bregman farthest-distance function and associated
farthest-point map and Klee sets. These objects were analyzed in
\cite{bwyy2} and are not treated here. In fact, under additional
assumptions, right and left notions may be related to each other via
duality. However, the duality approach was not powerful enough to
settle the question,  raised in \cite[Remark~7.3]{bwyy2},
\emph{whether or not every \fDK\ set is a singleton when $f$ does not
have full domain} as is the case when $D$ is, e.g., the
Kullback-Leibler divergence or the Itakura-Saito distance.
\emph{The first contribution of this paper is to settle this question
entirely}, for manifestations of $f$ that are even more general than those
considered in \cite{bwyy2}. In fact, in Theorem~\ref{t:super} we prove
that the answer is affirmative in the present setting.

Another related line of work concerns \emph{Chebyshev centers}.
Again, let us start by reviewing the classical situation in Euclidean
spaces. Let $S$ be a nonempty compact subset of $\RR^J$. The
\emph{Chebyshev center} is the center of the
smallest closed ball one can place in $\RR^J$ that entirely
captures the set $S$. The Chebyshev center exists and is unique,
and a classical result due to Garkavi and Klee
(see Corollary~\ref{c:G-K} below) provides a geometric characterization
of it. Unlike Klee sets, Chebyshev centers have already been
investigated in the context of Bregman distances ---
see, e.g., the work by Nielsen and Nock \cite{NielsenNock,NockNielsen}
(see Corollary~\ref{c:N-N} below)
and the references therein --- however;
it is assumed there that $S$ is \emph{finite}.
\emph{The second contribution of this paper is to extend the
classical work in Euclidean space by Garkavi and Klee and
the recent work by Nielsen and Nock on Chebyshev centers of finite
sets with respect to Bregman distances}. In Theorem~\ref{t:chebcent}, we
prove existence and uniqueness for Chebyshev centers of compact sets
with respect the Bregman distance, and we present a geometric
characterization of it.

The remainder of the paper is organized as follows.
In Section~\ref{s:aux} we collect and present several results
that will make the proofs of the main results more structured and
easier to follow. The main result in Section~\ref{s:main1} is
Theorem~\ref{t:super}, which states that every compact \fDK\ set
is indeed a singleton. In Section~\ref{s:main2}, we guarantee
existence and uniqueness of the $\fD{}\,$-Chebyshev center,
and we characterize it geometrically.
In Section~\ref{s:chebcent},
we illustrate our results with an example for three Bregman distances.

\section{Auxiliary Results}

\label{s:aux}

In this section, we collect several
results that will make the proofs of the main results
easier to follow. We start with two identities that
are straightforward consequences of \eqref{eq:D}.

\begin{lemma} \emph{(See \cite[Lemma~3.1]{CT}.)}
\label{l:3points}
Let $x$ be in $\RR^J$, and let $y$ and $z$ be in $U$.
Then
\begin{equation}
\label{e:3points}
D(x,z) - D(y,z) = D(x,y) + \scal{x-y}{\nabla f(y)-\nabla f(z)}.
\end{equation}
\end{lemma}

\begin{lemma} \emph{(See \cite[Remark~2.5]{BL}.)}
\label{l:4points}
Let $x_1$ and $x_2$ be in $\dom f$, and let $y_1$ and $y_2$ be in $U$.
Then
\begin{equation}
\label{e:4points}
\scal{x_1-x_2}{\nabla f(y_1)-\nabla f(y_2)} = D(x_2,y_1) + D(x_1,y_2)-D(x_1,y_1)-D(x_2,y_2).
\end{equation}
\end{lemma}

\begin{lemma}
\label{l:Dcont}
The Bregman distance $D$ is continuous on $U\times U$.
\end{lemma}
\begin{proof}
This follows from \cite[Theorem~10.1 and Corollary~25.5.1]{Rock70}.
\end{proof}

\begin{fact}[Rockafellar] \emph{(See \cite[Theorem~26.5]{Rock70}.)}
\label{f:Legendre}
The gradient operator $\nabla f$ is a continuous bijection between
$U$ and $\intdom f^*$, with continuous inverse $(\nabla f)^{-1} = \nabla f^*$.
Furthermore, $f^*$ is also a convex function of Legendre type.
\end{fact}

Recall that a function $g\colon\RR^J\to\RX$ is \emph{coercive} if
all its lower level sets are bounded; equivalently, if $\lim_{\|x\|\to\pinf} g(x)=\pinf$.
The following is thus clear.
\begin{equation}
\label{e:Moro:ii}
\text{If $g\colon\RR^J\to\RX$ is coercive and lower semicontinuous,
then $\argmin g\neq\varnothing$.}
\end{equation}
Here $\argmin g$ denotes the set of minimizers of $g$.

\begin{fact} \emph{(See \cite[Corollary~14.2.2]{Rock70}.)}
\label{f:MoRo}
Let $g\colon \RR^J\to\RX$ be convex, lower semicontinuous, and proper,
and let $x^*\in\RR^J$.
Then
$g(\,\cdot\,)-\scal{\,\cdot\,}{x^*}$ is coercive if and only if
$x^*\in\intdom g^*$.
\end{fact}

\begin{fact}[Ioffe-Tikhomirov]
\label{f:IT}
\emph{(See \cite[Theorem~2.4.18]{Zali02}.)}
Let $A$ be a compact Hausdorff space,
let $g_a\colon \RR^J\to\RX$ be convex for every $a\in A$,
and set $g := \sup_{a\in A} g_a$.
Assume that $(\forall x\in\RR^J)$
$A \to \RX\colon a \mapsto g_a(x)$ is upper semicontinuous
and that $x_0\in\dom g$ is a point such that
$(\forall a\in A)$ $g_a$ is continuous at $x_0$.
Then
\begin{equation}
\partial g(x_0) = \cconv \bigcup_{\menge{a\in A}{g(x_0)=g_a(x_0)}} \partial g_a(x_0).
\end{equation}
\end{fact}

\begin{lemma} \emph{(See \cite[Proposition~3.16]{Baus97}.)}
\label{l:leftBproj}
Suppose that $C$ is closed and convex, and let $y\in U\smallsetminus C$.
Then there exists a unique point $\bar{c} \in C$ such that
\begin{equation}
\big(\forall c\in C\big)\quad
\scal{c-\bar{c}}{\nabla f(y)-\nabla f(\bar{c})} \leq 0.
\end{equation}
\end{lemma}

\begin{lemma} \label{l:key}
Suppose that $C$ is compact, and
let $x\in U\smallsetminus (  (\nabla f^*)(\cconv(\nabla f(C))))$.
Then there exists $y\in (\nabla f^*)(\cconv(\nabla f(C)))\subset U$
such that
\begin{equation}
\big(\forall c\in C\big)\quad
D(x,c)\geq D(x,y)+D(y,c).
\end{equation}
\end{lemma}
\begin{proof}
Set $S := \nabla f(C)$ and $V := \intdom f^* = \nabla f(U)$.
Since $C$ is compact and $\nabla f$ is continuous (Fact~\ref{f:Legendre}),
the set $S$ is compact. Using \cite[Theorem~17.2]{Rock70}, we deduce that
$\conv S = \cconv S$ is a nonempty proper compact subset of $V$.
\opt{hhb}{{\color{blue}\marginpar{\color{blue}hhb}
\\$\blacktriangleright$
$V$ is an open set and $\cconv S$ is a closed bounded set.
The only sets in $\RR^J$ that are both open and closed
are $\varnothing$ and $\RR^J$.
$\blacktriangleleft$
}}
Using Fact~\ref{f:Legendre} again, we see that
\begin{equation}
 (\nabla f^*)(\cconv S)\;\;
 \text{is a proper compact subset of $U$}
\end{equation}
and that $x^* := \nabla f(x) \in V\smallsetminus (\cconv S)$.
Applying Lemma~\ref{l:leftBproj} (to $f^*$, $\cconv S$, and $x^*$), we obtain
a point $y^*\in \cconv S$ such that
\begin{equation}
\label{e:wandern1}
\big(\forall v \in \cconv S\big) \quad
\scal{v - y^*}{\nabla f^*(x^*)-\nabla f^*(y^*)} \leq 0.
\end{equation}
Now set $y := \nabla f^*(y^*)$. Then
\eqref{e:wandern1} yields
\begin{equation}
\big(\forall c\in C\big)\quad
\scal{\nabla f(y)-\nabla f(c)}{x-y} \geq 0.
\end{equation}
Combining this with Lemma~\ref{l:3points}, we estimate
\begin{equation}
\big(\forall c\in C\big)\quad
D(x,c)-D(y,c) = D(x,y) + \scal{\nabla f(y)-\nabla f(c)}{x-y}
\geq D(x,y),
\end{equation}
which completes the proof.
\end{proof}

Let $X$ and $Y$ be nonempty subsets of $\RR^J$ and let
$A\colon X\To Y$ be a set-valued operator, i.e.,
$(\forall x\in X)$ $Ax\subseteq Y$.
Denote the \emph{graph} of $A$ by
$\gr A := \menge{(x,y)\in X\times Y}{y\in Ax}$.
We say that $A$ is \emph{monotone from $X$ to $Y$}, if
\begin{equation}
\big(\forall (x,x^*)\in\gr A\big)
\big(\forall (y,y^*)\in\gr A\big)
\quad
\scal{x-y}{x^*-y^*} \geq 0.
\end{equation}
If $A$ is monotone from $X$ to $Y$ and every proper set-valued extension
from $X$ to $Y$  is not monotone, then $A$ is \emph{maximal monotone
from $X$ to $Y$}.
If $X=Y=\RR^J$, we will simply speak of monotone and maximal
monotone operators; this is the usual and well known setting.

We now present a variant of \cite[Example~12.7]{Rock98}, which
is a sufficient condition for maximal monotonicity.

\begin{proposition}\label{p:maximal}
Let $O$ be a nonempty open subset of $\RR^J$,
let $Y$ be a subset of $\RR^J$, and let
$A\colon O\to Y$ be monotone and continuous.
Then $A$ is maximal monotone from $O$ to $Y$.
\end{proposition}
\begin{proof}
Suppose that $(\bar{x},\bar{y})\in O\times Y$ satisfies
\begin{equation}
\big(\forall x\in O\big)\quad \scal{\bar{x}-x}{\bar{y}-A{x}}\geq 0,
\end{equation}
and denote the closed unit ball in $\RR^J$ by $B$.
Then for all sufficiently small $\epsilon>0$, we have
$\bar{x}+\epsilon B \subseteq U$ and hence
$(\forall b\in B)$
$\scal{\bar{x}-(\bar{x}+\epsilon b)}{\bar{y}-A(\bar{x}+\epsilon b)}\geq 0$
and so $\scal{b}{\bar{y}-A(\bar{x}+\epsilon b)}\leq 0$.
Letting $\epsilon\to 0^+$ for fixed but arbitrary $b\in B$, and using
continuity of $A$ at $\bar{x}$, we deduce that
$\scal{b}{\bar{y}-A\bar{x}}\leq 0$. Supremizing this last inequality over
$b\in B$, we obtain $\|\bar{y}-A\bar{x}\|=0$.
Hence $(\bar{x},\bar{y}) = (\bar{x},A\bar{x})\in\gr A$,
as required.
\end{proof}

Our first result reveals a monotonicity property of $\ffproj{C}$.
(See also \cite{WesSch} and \cite[Proposition 7.1]{bwyy2},
and \cite{bwyy1}, where we discuss Chebyshev sets instead of Klee sets.)

\begin{proposition}\label{p:monotone}
The set-valued operator $-\nabla f \circ \ffproj{C}\colon\RR^J\To\RR^J$ is monotone.
\end{proposition}

\begin{proof}
Assume that $(x,x^*)$ and $(y,y^*)$ lie in $\gr\ffproj{C}$.
It follows from \eqref{e:montag:b} and Lemma~\ref{l:4points}
(applied to $x_1=x$, $x_2=y$, $y_1 = y^*$, and $y_2=x^*$) that
\begin{align}
0 &\leq \big(D(x,x^*)-D(x,y^*)\big) + \big(D(y,y^*)-D(y,x^*)\big)\\
&= \scal{x-y}{\nabla f(y^*)-\nabla f(x^*)}\notag\\
&= \scal{x-y}{(-\nabla f)(x^*)-(-\nabla f)(y^*)},\notag
\end{align}
as required.
\end{proof}

\begin{proposition}\label{p:cont}
Suppose that $C$ is closed, and that $((x_n,y_n))_\nnn$
is a sequence in $(\gr \ffproj{C}) \cap (U\times\RR^J)$
such that $(x_n,y_n)\to (x,y)\in U\times\RR^J$.
Then $(x,y)\in\gr \ffproj{C}$.
\end{proposition}
\begin{proof}
Since $\ran \ffproj{C} \subseteq C$, the sequence $(y_n)_\nnn$
lies in $C$ and it satisfies
$(\forall\nnn)$ $D(x_{n},y_{n})=\ffD{C}(x_{n})$.
Because $C$ is closed, $y\in C \subseteq U$.
By Lemma~\ref{l:Dcont}, $D$ is continuous on $U\times U$.
In view of \eqref{e:090217:a}, we deduce altogether
\begin{equation}
\ffD{C}(x)\leq\varliminf_{\nnn} \ffD{C}(x_{n})=
\varliminf_{\nnn} D(x_{n},y_{n})
=D(x,y)\leq \ffD{C}(x).
\end{equation}
Therefore, $\ffD{C}(x) = D(x,y)$, i.e., $y\in\ffproj{C}(x)$.
\end{proof}

\begin{proposition}
\label{p:090217}
Suppose that $\overline{C}\subseteq U$.
Then $\gr \ffproj{C} \subseteq \gr \ffproj{\overline{C}}$.
\end{proposition}
\begin{proof}
Take $(x,y)\in\gr  \ffproj{C}$. Then $y\in C \subseteq \overline{C}$
and $(\forall c\in C)$ $D(x,y)\geq D(x,c)$.
Since $\overline{C}\subseteq U$ and $D(x,\cdot)$ is continuous on $U$
(Lemma~\ref{l:Dcont}), it follows that
$(\forall \bar{c}\in \overline{C})$ $D(x,y)\geq D(x,\bar{c})$.
Thus, $y\in\ffproj{\overline{C}}(x)$.
\end{proof}

\begin{remark}
\label{r:090317}
Assume that $\bar{c} \in \overline{C} \cap \bd U$.
In view of \eqref{e:CsubU}, there
exists a sequence $(c_n)_\nnn$ in $C \subseteq U$ such
that $c_n\to \bar{c}$; hence,
by \cite[Theorem~3.8.(i)]{Baus97}, $(\forall x\in U)$ $D(x,c_n)\to\pinf$.
Therefore, the assumption that $\overline{C}$ be a subset of $U$ is very natural
in Proposition~\ref{p:090217} and elsewhere in this paper.
\end{remark}

\begin{proposition}
\label{p:090217n}
Suppose that $(\forall x\in \dom f)$ $D(x,\,\cdot\,)$ is convex on $U$.
Then $\gr\ffproj{C} \subseteq\gr\ffproj{\conv{C}}$.
\end{proposition}
\begin{proof}
Take $(x,y)\in\gr  \ffproj{C}$. Then $x\in\dom f$,
$y\in C \subseteq \conv{C}$,
and $(\forall c\in C)$ $D(x,c)\leq D(x,y)$.
Now let $z\in\conv C$, say $z=\sum_{i=1}^n\lambda_ic_i$,
where each $\lambda_i\in [0,1]$, each $c_i\in C$, and
$\sum_{i=1}^n \lambda_i = 1$.
Then $D(x,z) \leq \sum_{i=1}^{n}\lambda_iD(x,c_i)
\leq \sum_{i=1}^{n}\lambda_iD(x,y) = D(x,y)$ and therefore
$y\in\ffproj{\conv C}(x)$.
\end{proof}

\begin{proposition}
\label{p:090217nn}
Suppose that $(\forall x\in \dom f)$ $D(x,\,\cdot\,)$ is
strictly convex on $U$.
Then $\gr\ffproj{C} = \gr\ffproj{\conv{C}}$.
\end{proposition}
\begin{proof}
In view of Proposition~\ref{p:090217n},
we only need to show that
$\gr\ffproj{\conv C} \subseteq\gr\ffproj{{C}}$.
To this end, let $(x,y)\in\gr  \ffproj{\conv C}$.
Then $x\in \dom f$,
$y\in \conv{C}$,
and $(\forall s\in \conv C)$ $D(x,s)\leq D(x,y)$.
In particular,
$(\forall c\in C)$ $D(x,c)\leq D(x,y)$.
The proof is complete as soon as we have verified that
$y\in C$.
Assume to the contrary that $y\notin C$.
Then $y=\sum_{i=1}^{n}\lambda_i c_i$,
where $n\geq 2$, each $\lambda_i>0$, each $c_i\in C$,
and where the $c_i$ are pairwise distinct and $\sum_{i=1}^{n}\lambda_i=1$.
But then $D(x,y) < \sum_{i=1}^{n}\lambda_i D(x,c_i)
\leq\sum_{i=1}^{n}\lambda_i D(x,y) = D(x,y)$, which is absurd.
\end{proof}

The next result shows that when $D$ is separately
convex (see \cite{BB01} for a systematic discussion of separate and joint
convexity of $D$), then the
farthest-point distance is ``blind'' to the convex hull.
\begin{proposition}
Suppose that $(\forall x\in\dom f)$ $D(x,\,\cdot\,)$ is convex.
Then $\ffD{\conv C} = \ffD{C}$.
\end{proposition}
\begin{proof}
This follows from \cite[Theorem~32.2]{Rock70}.
\end{proof}

\section{Klee Sets are Singletons}

\label{s:main1}

The following result will be critical in the proof of our
first main result (Theorem~\ref{t:super}).

\begin{theorem}
\label{t:main}
Suppose that $C$ is compact.
Then $\argmin\ffD{C}$ is a nonempty subset of $U$.
\end{theorem}
\begin{proof}
By \eqref{e:090322n},
$\dom\ffD{C} = \dom f$.
Since $C\subset U$, it follows from Fact~\ref{f:Legendre}
 that $\nabla f(C) \subset \nabla f(U) = \intdom f^*$.
In view of Fact~\ref{f:MoRo}, we deduce that
\begin{equation}
\big(\forall c\in C\big) \quad
f(\,\cdot\,) - \scal{\,\cdot\,}{\nabla f(c)}\;\;\text{is coercive.}
\end{equation}
Since $(\forall c\in C)$ $D(\,\cdot\,,c) = ( f(\,\cdot\,) - \scal{\,\cdot\,}{\nabla f(c)} )
+ (\scal{c}{\nabla f(c)} - f(c))$, it follows that
\begin{equation}
\big(\forall c\in C\big) \quad
\;D(\,\cdot\,,c)\;\;\text{is coercive.}
\end{equation}
In turn, this implies that
\begin{equation}
\label{e:090216:a}
\ffD{C}(\,\cdot\,) = \sup_{c\in C}D(\,\cdot\,,c)\;\;
\text{is coercive, convex, lower semicontinuous, and proper.}
\end{equation}
\opt{hhb}{{\color{blue}\marginpar{\color{blue}hhb}
\\$\blacktriangleright$
supremum of coercive is coercive; supremum of convex is convex;
supremum of lsc is lsc; proper because $\dom \ffD{C} = \dom f$.
$\blacktriangleleft$
}}
In view of \eqref{e:090216:a} and \eqref{e:Moro:ii},
$\argmin\ffD{C}\neq\varnothing$. Let
\begin{equation} \label{e:090216:g}
x_0\in\argmin\ffD{C}.
\end{equation}
It suffices to show that
\begin{equation}
\label{e:090216:gg}
x_0\in U.
\end{equation}
Assume to the contrary that $x_0\notin U$.
In view of \eqref{e:090216:f} and \eqref{e:090216:g},
$x_0 \in (\dom f \smallsetminus U)\subseteq \bd\dom f$.
Now fix an arbitrary point $x_1\in U$ and set
\begin{equation} \label{e:acclemmma}
\big(\forall \epsilon\in\left]0,1\right[\big)\quad x_\epsilon :=
(1-\epsilon)x_0+\epsilon x_1.
\end{equation}
By \cite[Theorem~6.1]{Rock70},
$(\forall \epsilon\in\left]0,1\right])$, $x_\epsilon \in U$.
Set $S := \nabla f(C)$.
As already observed in the proof of Lemma~\ref{l:key},
$\conv S = \cconv S$ is a proper compact subset of $\intdom f^*$.
Thus, there exists $\bar{\epsilon}\in\left]0,1\right]$ such that
$(\forall \epsilon\in\left]0,\bar{\epsilon}\right])$
$x_\epsilon \in U\smallsetminus(\nabla f^*)(\cconv S)$.
\opt{hhb}{{\color{blue}\marginpar{\color{blue}hhb}
\\$\blacktriangleright$ For otherwise, we would find a sequence $(\epsilon_n)$ converging
to $0^+$ such that $x_{\epsilon_n} \in (\nabla f^*)(\cconv S)$.
Then $x^*_{n} := \nabla f(x_{\epsilon_n}) \in \cconv S$.
After passing to a subsequence if necessary, we assume
that $x_n^* \to x^*\in\cconv S\subset \intdom f^*$.
But then $x_{\varepsilon_n} \to \nabla f^*(x^*)\in U$.
On the other hand, by \eqref{e:acclemmma}, $x_{\epsilon_n} \to x_0
\notin U$ --- contradiction!!
$\blacktriangleleft$
}}
Lemma~\ref{l:key} now yields
\begin{equation} \label{e:wandern2}
\big(\forall\epsilon\in\left]0,\bar{\epsilon}\right]\big)
\big(\exists y_\epsilon\in (\nabla f^*)(\cconv S)\big)
\big(\forall c\in C\big)
\quad D(x_\epsilon,c) \geq D(x_\epsilon,y_\epsilon)
+D(y_\epsilon,c).
\end{equation}
On the one hand, while $f$ is not necessarily continuous at $x_0$, it is at
least \emph{continuous along the line segment} $[x_0,x_1]$
(see \cite[Theorem~7.5]{Rock70}); consequently,
\begin{equation}
\label{e:090216:l1}
\lim_{\epsilon\to 0^+} f(x_\epsilon) = f(x_0).
\end{equation}
On the other hand, the net
$(y_\epsilon)_{\epsilon\in\left]0,\bar{\epsilon}\right]}$ lies in
$\nabla f^*(\cconv S)$, which is a compact set.
After passing to a subnet and relabeling if necessary, we assume
that there exists a point $y_0\in\RR^J$ such that
\begin{equation}
\label{e:090216:l2}
\lim_{\epsilon\to 0^+} y_\epsilon = y_0 \in \nabla f^*(\cconv S)\subset U.
\end{equation}
Combining \eqref{e:090216:l1} and \eqref{e:090216:l2},
invoking Lemma~\ref{l:Dcont},
and taking the limit in \eqref{e:wandern2}, we obtain altogether that
\begin{equation} \label{e:090216:s}
\big(\forall c\in C\big)\quad
D(x_0,c) \geq D(x_0,y_0) + D(y_0,c).
\end{equation}
Since $x_0\in\bd\dom f$ and $y_0\in\intdom f = U$, \eqref{e:misnomer}
results in $D(x_0,y_0)>0$.
Supremizing \eqref{e:090216:s} over $c\in C$, we deduce
that
\begin{equation}
\ffD{C}(x_0) \geq D(x_0,y_0) + \ffD{C}(y_0) > \ffD{C}(y_0),
\end{equation}
which contradicts \eqref{e:090216:g}.
Therefore, we have verified \eqref{e:090216:gg},
and the proof is complete.
\end{proof}

\begin{theorem}[every \fDK\ set is a singleton]
\label{t:super}
Suppose that $C$ is compact and
that $C$ is \fDK. Then $C$ is a singleton.
\end{theorem}
\begin{proof}
Recall that
\begin{equation}
\ffD{C}(\,\cdot\,) = \sup_{c\in C} D(\,\cdot\,,c) =
\sup_{c\in C} \Big(\big(f(\,\cdot\,) - \scal{\,\cdot\,}{\nabla f(c)}\big) +
\big(\scal{c}{\nabla f(c)} - f(c)\big)\Big).
\end{equation}
Because $C$ is \fDK,
if $x\in U$, then $\ffproj{C}x$ is the unique point in $C$ such that
$\ffD{C}(x) = D(x,\ffproj{C}x)$ and
$(\forall c\in C\smallsetminus\{\ffproj{C}x\})$
$\ffD{C}(x)>D(x,c)$.
In view of Theorem~\ref{t:main}, we take $x_0\in\argmin\ffD{C}\subset U$.
Using the Fact~\ref{f:IT}, we obtain
\begin{equation}
0\in \partial \ffD{C}(x_0) =  \nabla f(x_0)-\nabla f(\ffproj{C}x_0).
\end{equation}
\opt{hhb}{{\color{blue}\marginpar{\color{blue}hhb}
\\$\blacktriangleright$ We apply Fact~\ref{f:IT} with $A = C$, which is compact by hypothesis.
Furthermore $g_a = g_c = D(\cdot,c)$.
Then $g= \ffD{C}$ and $\dom g = \dom \ffD{C} = \dom f$.
The point $x_0$ belongs to $U=\intdom g$.
Since $g$ is convex, $g$ is continuous at $x_0$,
and the same is true for every $g_a = D(\cdot,a)$.
For every $x\in\RR^J$, the function
$C\to\RX\colon c\mapsto g_c(x) = D(x,c)
= f(x) - f(c)-\scal{x-c}{\nabla f(c)}$
is certainly upper semicontinuous in $C$;
indeed, this function is either $\equiv\pinf$ if $x\notin\dom f$;
or it is clearly continuous if $x\in\dom f$.
$\blacktriangleleft$
}}

Hence $\nabla f(x_0) = \nabla f(\ffproj{C}(x_0))$ and thus
$x_0 = \ffproj{C}(x_0)$.
Therefore, $C = \{x_0\}$.
\end{proof}

\begin{corollary}[Klee]
Suppose that $C$ is compact Klee set with respect to the
Euclidean distance. Then $C$ is a singleton.
\end{corollary}
\begin{proof}
(See also \cite{Klee61}.)
This follows from Theorem~\ref{t:super} when $f = \thalb\|\cdot\|^2$.
\end{proof}

We conclude this section with two results
concerning \fDK\ sets that are not assumed to be compact.
When considering classical Klee sets,
a standard assumption is \emph{closedness}.
The next result illustrates this assumption
in the present Bregman distance setting.

\begin{proposition}
Suppose that $\overline{C}$ is a compact subset of $U$,
and that $U \subseteq \dom \ffproj{C}$.
Then
$\overline{C}$ is \fDK\
if and only if $C$ is \fDK\ and $\ffproj{C}$ is continuous on $U$.
\end{proposition}
\opt{hhb}{{\color{blue}\marginpar{\color{blue}hhb}
\medskip
$\blacktriangleright$
This leaves open the door for a bounded nonclosed \fDK\ set $C$
such that $\ffproj{C}$ is not continuous??
$\blacktriangleleft$

\medskip
}}
\begin{proof}
``$\Rightarrow$'':
Since $\overline{C}$ is compact, Theorem~\ref{t:super}
implies that $\overline{C}$ is a singleton, say $\overline{C}=\{y\}$.
But then $C=\{y\}=\overline{C}$ is also \fDK, and
$\ffproj{C}|_U \equiv\{y\}$ is clearly continuous on $U$.

``$\Leftarrow$'': Proposition~\ref{p:monotone} implies that
both
$-\nabla f\circ \ffproj{C}\big|_U $ and
$-\nabla f\circ \ffproj{\overline{C}}\big|_U$ are monotone from
$U$ to $\RR^J$.
Furthermore, since $\ffproj{C}$ is continuous on $U$, so is
$-\nabla f\circ \ffproj{C}\big|_U$.
Thus, by Proposition~\ref{p:maximal},
$-\nabla f\circ \ffproj{C}\big|_U$ is maximal monotone from $U$
to $\RR^J$.
On the other hand,
Proposition~\ref{p:090217} implies that
\begin{equation}
\gr\big(-\nabla f\circ \ffproj{C}\big|_U\big) \subseteq
\gr\big(-\nabla f\circ \ffproj{\overline{C}}\big|_U\big).
\end{equation}
Altogether, $-\nabla f\circ \ffproj{C}\big|_U = -\nabla f\circ
\ffproj{\overline{C}}\big|_U$, which yields
$\ffproj{C}\big|_U = \ffproj{\overline{C}}\big|_U$.
Since $\ffproj{C}\big|_U$ is single-valued, so is
$\ffproj{\overline{C}}\big|_U$. Therefore, $\overline{C}$
is \fDK.
\end{proof}

If the underlying Bregman distance $D$ is strictly convex in the second
variable, then we obtain the following result.

\begin{proposition}
Suppose that $(\forall x\in U)$ $D(x,\,\cdot\,)$ is strictly convex on $U$.
Then $\conv C$ is \fDK\ if and only if $C$ is \fDK.
\end{proposition}
\begin{proof}
This is an immediate consequence of Proposition~\ref{p:090217nn}.
\end{proof}

\section{Characterization of Chebyshev Centers}

\label{s:main2}

The proof of our second main result (Theorem~\ref{t:chebcent})
relies upon the next two results.

\begin{proposition}
\label{p:strictcon}
Suppose that $C$ is compact.
Then $\ffD{C}$ is proper, lower semicontinuous, and convex,
with $\dom \ffD{C} = \dom f = \dom\ffproj{C}$. Furthermore,
$\ffD{C}$ is strictly convex on $\dom \partial \ffD{C} = \intdom f = U$.
\end{proposition}
\begin{proof}
We observed already (see \eqref{e:090217:a} and \eqref{e:090322n})
that $\ffD{C}$ is convex and lower semicontinuous,
and that
$\dom \ffD{C} = \dom f = \dom\ffproj{C}$.
Hence $\ffD{C}$ is proper.
Now set
\begin{equation}
g\colon \RR^J\to\RX\colon
x \mapsto \max_{c\in C}\big( \scal{c-x}{\nabla f(c)} - f(c)\big),
\end{equation}
and note that $g$ is convex with $\dom g=\RR^J = \intdom\partial g$
 (see \cite[Theorem~23.4]{Rock70}).
Furthermore,
\begin{equation}
\ffD{C} = f + g.
\end{equation}
By the subdifferential sum rule (see \cite[Theorem~23.8]{Rock70}), we
have $\partial \ffD{C} = \partial f + \partial g$ and hence
$\dom \partial \ffD{C} = \dom(\partial f) \cap \dom(\partial g)
=\dom(\partial f) \cap \RR^J = \dom\partial f$.
On the other hand, since $f$ is a Legendre function, it follows
from \cite[Theorem~26.1]{Rock70} that
$\dom\partial f = \intdom f$.
Altogether, $\dom\partial \ffD{C}=\intdom f = U$.
Using once more the assumption that $f$ is a Legendre function,
we have that $f$ is strictly convex on $\intdom f = U$, and therefore
so is $\ffD{C} = f+g$.
\end{proof}

Recall that for a proper convex function $g:\RR^{J}\to\RX$,
the \emph{directional derivative}
of $g$ at $x\in\dom g$ in the direction $h\in\RR^{J}$ is defined by
\begin{equation}
g'(x;h)=\lim_{t\to 0^+}\frac{g(x+th)-g(x)}{t}.
\end{equation}

\begin{theorem}[directional derivative]
\label{t:dirdev}
Suppose that $C$ is compact,
let $x\in\dom f$, and let $h\in \RR^J$.
Then
\begin{equation}\label{directional}
\ffD{C}'(x;h)=
\sup\menge{f'(x;h)-\scal{{h}}{\nabla f(y)}}{y\in \ffproj{C}(x)}.
\end{equation}
If $x\notin U$ and $x+h\in U$, then $\ffD{C}'(x;h)=\minf$.
\end{theorem}
\begin{proof}
Recall that $\dom\ffD{C} = \dom f = \dom\ffproj{C}$
by Proposition~\ref{p:strictcon}, so let $y\in\ffproj{C}(x)$.
Then
\begin{equation}
\big(\forall t>0\big)\quad
\ffD{C}(x+th)\geq D(x+th,y) = f(x+th)-f(y)-\scal{x+th-y}{\nabla f(y)}
\end{equation}
and
\begin{equation}
\ffD{C}(x)=D(x,y) = f(x)-f(y)-\scal{x-y}{\nabla f(y)}.
\end{equation}
Hence, $(\forall t>0)$
$\ffD{C}(x+th)-\ffD{C}(x)\geq f(x+th)-f(x)-\scal{th}{\nabla
f(y)}$.
Dividing by $t$ and taking the infimum over $t>0$ yields
\begin{equation}
\ffD{C}'(x;h) \geq f'(x;h) - \scal{h}{\nabla f(y)}.
\end{equation}
Supremizing over $y\in\ffproj{C}(x)$ yields
\begin{equation}
\label{e:catscan:1}
\ffD{C}'(x;h) \geq \sup\menge{f'(x;h) - \scal{h}{\nabla
f(y)}}{y\in\ffproj{C}(x)}.
\end{equation}
If $[x,x+h]\cap\dom f=\{x\}$, then $f'(x;h)=\pinf$; hence,
\eqref{directional} follows from \eqref{e:catscan:1}.
Thus, we assume that $[x,x+h] \cap \dom f$
contains a nontrivial line segment.
Let $(t_n)_\nnn$ be a sequence in $\left]0,1\right[$ such that
$t_n\to 0^+$ and $(x+t_nh)_\nnn$ lies in $\dom f$.
Furthermore, for every $\nnn$, let $c_n \in \ffproj{C}(x+t_nh)$.
After passing to a subsequence and relabeling if necessary,
we also assume that $c_n \to \bar{c}\in C$.
Then, for every $\nnn$,
\begin{equation}
\label{e:catscan:3}
\ffD{C}(x+t_nh) = D(x+t_nh,c_n) = f(x+t_nh)-f(c_n)-\scal{x+t_nh-c_n}{\nabla
f(c_n)}
\end{equation}
and
$\ffD{C}(x) \geq D(x,c_n) = f(x)-f(c_n)-\scal{x-c_n}{\nabla f(c_n)}$;
consequently,
\begin{equation} \label{e:catscan:2}
\frac{\ffD{C}(x+t_{n}h)-\ffD{C}(x)}{t_{n}}\leq
\frac{f(x+t_{n}h)-f(x)}{t_{n}}-\scal{h}{\nabla f(c_{n})}.
\end{equation}
Letting $n\to\pinf$ in \eqref{e:catscan:2}, we deduce that
\begin{equation}
\label{e:catscan:4}
\ffD{C}'(x;h) \leq f'(x;h) - \scal{h}{\nabla f(\bar{c})}.
\end{equation}
On the other hand, using line segment continuity of $f$
and $\ffD{C}$ at $x$ (see \cite[Corollary~7.5.1]{Rock70}),
and continuity of both $f$ and $\nabla f$ on $U$, we see
that letting $n\to\pinf$ in \eqref{e:catscan:3} yields
$\ffD{C}(x) = D(x,\bar{c})$. Hence $\bar{c}\in\ffproj{C}(x)$.
It thus follows from \eqref{e:catscan:4} that
$\ffD{C}'(x;h) \leq
\sup\menge{f'(x;h)-\scal{{h}}{\nabla f(y)}}{y\in \ffproj{C}(x)}$.
Combining this with  \eqref{e:catscan:1}, we deduce that
\eqref{directional} holds.
The ``If'' statement follows from \eqref{directional} and
\cite[Theorem~23.3]{Rock70}.
\end{proof}

\begin{theorem}[subdifferential]
\label{t:subdiff}
Suppose that $C$ is compact, and let $x\in U$.
Then
\begin{equation}
\partial \ffD{C}(x) = \nabla f(x) - \conv\nabla f\big(\ffproj{C}(x)\big).
\end{equation}
\end{theorem}
\begin{proof}
By Theorem~\ref{t:dirdev} and
\cite[Theorem~23.4]{Rock70}, $\ffD{C}'(x;\cdot)$ is the support function of
both $\nabla f(x) - \nabla f(\ffproj{C}(x))$ and $\partial \ffD{C}(x)$.
Therefore, the latter set (which is closed and convex already) is
the closed convex hull of the former set.
Since $\ffproj{C}(x)$ is a compact subset of $U$ by Proposition~\ref{p:cont},
it follows from the continuity of $\nabla f$ on $U$ and from
\cite[Theorem~17.2]{Rock70} that
$\cconv \nabla f(\ffproj{C}(x)) = \conv \nabla f(\ffproj{C}(x))$.
This completes the proof.
\end{proof}

\begin{theorem}[uniqueness and characterization of the
$\fD{}\,$-Chebyshev center]
\label{t:chebcent}
Suppose that $C$ is compact.
Then  $\ffD{C}$ has a unique minimizer $x\in\dom f$, called the
\emph{$\fD{}\,$-Chebyshev center} of $C$, and characterized by
\begin{equation}\label{e:chebcent}
x\in \nabla f^*\Big(\conv \nabla f\big(\ffproj{C}(x)\big)\Big).
\end{equation}
\end{theorem}
\begin{proof}
Theorem~\ref{t:main} states that $\argmin\ffD{C}$ is a nonempty
subset of $U$. In view of the strict convexity of $\ffD{C}$ on $U$
 (Proposition~\ref{p:strictcon}),
$\argmin\ffD{C}$ is a singleton, say $\{x\}$.
By Theorem~\ref{t:subdiff},
$0\in\partial\ffD{C}(x) = \nabla f(x) - \conv\nabla f(\ffproj{C}(x))$ and
thus $\nabla f(x)\in\conv\nabla f(\ffproj{C}(x))$.
Now apply Fact~\ref{f:Legendre}.
\end{proof}

\begin{corollary}[Garkavi-Klee]
\label{c:G-K}
\emph{(See \cite{Garkavi} and also \cite{Klee60}.)}
Suppose that $C$ is compact and that $x\in\RR^J$.
Then $x$ is the Chebyshev center of $C$ with respect to
the Euclidean distance if and only if
\begin{equation}
x\in\conv\ffproj{C}(x).
\end{equation}
\end{corollary}

\begin{corollary}[Nock-Nielsen]
\label{c:N-N}
\emph{(See \cite{NockNielsen} and also \cite{NielsenNock}.)}
Suppose that $C$ is finite.
Then the $\fD{}\,$-Chebyshev center of $C$
is the unique point $x\in U$ characterized by
\begin{equation}
x\in\nabla f^*\Big(\conv \nabla f\big(\ffproj{C}(x)\big)\Big).
\end{equation}
\end{corollary}

\begin{corollary}
\label{c:chebmultval}
Suppose that $C$ is compact and that it contains at least $2$ points,
and let $x\in U$ be the $\fD{}\,$-Chebyshev center of $C$.
Then $\ffproj{C}(x)$ must contain at least $2$ points.
\end{corollary}
\begin{proof}
Suppose to the contrary that
$\ffproj{C}(x)$ is a singleton. Then
\eqref{e:chebcent} implies that $\ffproj{C}(x)=\{x\}$, i.e.,
that $x$ is its own farthest point in $C$.
In view of \eqref{e:misnomer} and the assumption that $C$ contains a point
different from $x$, this is absurd.
\end{proof}

\section{Constructing and Visualizing Chebyshev Centers}
\label{s:chebcent}

We work in the Euclidean plane, i.e., we assume that $J=2$,
and we let $D$ be the halved Euclidean
distance squared, the Kullback-Leibler divergence, or the Itakura-Saito
distance (see Example~\ref{ex:distances}).
Set
\begin{equation}
\bc_0 = (1,a)
\quad\text{and}\quad
\bc_1 = (a,1),
\quad
\text{where}
\quad
a\in\left]1,\pinf\right[,
\end{equation}
and
\begin{equation}
\big(\forall \lambda\in\RR\big)\quad
\bc_\lambda = (1-\lambda)\bc_0+\lambda\bc_1.
\end{equation}
Furthermore, we assume that
\begin{align}
C &=  \conv\{\bc_0,\bc_1\}
= \menge{\bc_\lambda}{\lambda\in[0,1]}
=
\menge{\big(1-\lambda+\lambda a,(1-\lambda)a+\lambda\big)}{\lambda\in[0,1]}
\\
&=
\menge{\left((a-1)\lambda+1,(1-a)\lambda+a\right)}{\lambda\in[0,1]}.\notag
\end{align}
Note that $C\subset \RR^2_{++}\subseteq U$, and that $C$ is compact
and convex.
In view of Theorem~\ref{t:chebcent},
the $\fD{}\,$-Chebyshev center $\bz$ of $C$
is characterized by
\begin{equation}
\label{e:zchar}
\bz \in\nabla f^*\Big( \conv\nabla f\big(\ffproj{C}(\bz)\big)\Big).
\end{equation}
Our aim in this section is to determine $\bz$ and related objects,
and to visualize them.
It will be convenient to set
\begin{equation}
\Delta = \menge{(x,x)}{x\in\RR}.
\end{equation}

\begin{proposition}
\label{p:ondiag}
$\bz\in\Delta$.
\end{proposition}
\begin{proof}
For $\bx = (x_1,x_2)\in\RR^2$, set $\bx^\intercal = (x_2,x_1)$.
Observe that for the choices of $D$ considered in this section,
$(\forall \bx\in \RR^2)(\forall \by\in\RR^2)$ $D(\bx,\by) =
D(\bx^\intercal,\by^\intercal)$ and that $C^\intercal =
\menge{\bc^\intercal}{\bc\in C} = C$.
Thus,
$(\forall\bx\in\RR^2)$ $\ffD{C}(\bx) = \ffD{C}(\bx^\intercal)$.
Since $\bz$ is the \emph{unique} minimizer of $\ffD{C}$, we must have
that $\bz=\bz^\intercal$, i.e., that $\bz\in\Delta$.
\end{proof}

\begin{example}[halved Euclidean distance squared]
\label{ex:halfeuclid}
Suppose $D$ is as in Example~\ref{ex:distances}\ref{ex:distances:E}, and
let $\bx=(x_1,x_2)\in\RR^2$.
Then
\begin{align}\label{e:euclid_farpt}
\ffproj C(\bx) & =\begin{cases}
\{\bc_0\},&\text{if $x_2<x_1$;}\\
\{\bc_1\},&\text{if $x_1<x_2$;}\\
\{\bc_0,\bc_1\},&\text{if $x_1=x_2$,}
\end{cases}
\end{align}
and
$\bz=\bc_{1/2} = \big(\thalb(1+a),\thalb(1+a)\big)$.
\end{example}
\begin{proof}
Set
\begin{equation}
d_\bx\colon\RR\to\RP\colon \lambda\mapsto D(\bx,c_\lambda).
\end{equation}
Then for every $\lambda\in\RR$, we have
\begin{equation}
d_\bx(\lambda) =
(a-1)^2\lambda^2+(1-a)(x_1-x_2+a-1)\lambda+\frac{(x_1-1)^2+(x_2-a)^2}{2},
\end{equation}
\begin{equation}
d'_\bx(\lambda) = (x_1-x_2-1+a)(1-a)+2\lambda(1-a)^2
\;\text{and}\;
d''_\bx(\lambda) = 2(a-1)^2.
\end{equation}
Hence $\ffproj{C}(\bx)\subseteq\{\bc_0,\bc_1\}$.
Since
$d_\bx(0)-d_\bx(1) = (1-a)(x_2-x_1)$, we obtain \eqref{e:euclid_farpt}.
Furthermore, since $C$ is convex and $\bc_{1/2}\in\Delta$, we have
$\bc_{1/2} \in C = \conv\{\bc_0,\bc_1\} = \conv\ffproj{C}(\bc_{1/2})$.
Therefore, the characterization \eqref{e:zchar} of $\bz$ yields
$\bz=\bc_{1/2}$.
(Alternatively, one may verify that $\bc_{1/2}$ is the unique minimizer
of the function
$\Delta\to\RP\colon(x,x)\mapsto d_{(x,x)}(0) = \ffD{C}(x,x)$.)
\end{proof}

\begin{example}[Kullback-Leibler divergence]
\label{ex:kullback}
Suppose $D$ is as in Example~\ref{ex:distances}\ref{ex:distances:KL},
and let $\bx = (x_1,x_2)\in U$.
Then
\begin{align}\label{e:KL_farpt}
\ffproj C(\bx) & =\begin{cases}
\{\bc_0\},&\text{if $x_2<x_1$;}\\
\{\bc_1\},&\text{if $x_1<x_2$;}\\
\{\bc_0,\bc_1\},&\text{if $x_1=x_2$,}
\end{cases}
\end{align}
and
$\bz=\big(\sqrt{a},\sqrt{a}\,\big)$.
\end{example}
\begin{proof}
Set
\begin{equation}
d_\bx\colon\RR\to\left[0,\pinf\right]\colon\lambda\mapsto D(\bx,\bc_\lambda).
\end{equation}
Then $\dom d_\bx = \menge{\lambda\in\RR}{\bc_\lambda\in U} =
\left]-1/(a-1), a/(a-1)\right[ \supset [0,1]$.
For every $\lambda\in\dom d_\bx$, we have
\begin{equation}
d_\bx(\lambda) =
-x_1\ln\left(\frac{(a-1)\lambda+1}{x_1}\right)-x_1+1-x_2\ln\left(\frac{(1-a)\lambda+a}{x_2}\right)-x_2+a,
\end{equation}
\begin{equation}
d'_\bx(\lambda) =
-\frac{x_1(a-1)}{(a-1)\lambda+1}-\frac{x_2(1-a)}{(1-a)\lambda+a},
\end{equation}
and
\begin{equation}
d''_\bx(\lambda)
=\frac{x_1(a-1)^2}{\big((a-1)\lambda+1\big)^2}+
\frac{x_2(1-a)^2}{\big((1-a)\lambda+a\big)^2}>0.
\end{equation}
Thus, $d_\bx$ has no local maximizers in $\dom d_\bx$ and therefore
$\ffproj{C}(\bx)\subseteq\{\bc_0,\bc_1\}$.
Because of
\begin{equation}
D(\bx,\bc_0)-D(\bx,\bc_1) = d_\bx(0)-d_\bx(1) = (x_1-x_2)\ln(a),
\end{equation}
we see that \eqref{e:KL_farpt} must hold.
Finally, \eqref{e:KL_farpt} implies that
\begin{align}
\big(\sqrt{a},\sqrt{a}\,\big) &=
\bigg(\exp\Big(\thalb\big(0+\ln(a)\big)\Big),
\exp\Big(\thalb\big(\ln(a)+0\big)\Big)\bigg)\\
&= (\exp\times\exp)
\Big(\thalb\big(\ln(1),\ln(a)\big) + \thalb\big(\ln(a),\ln(1)\big)\Big)\notag\\
&= \nabla f^*
\big(\thalb\nabla f(\bc_0) + \thalb\nabla f(\bc_1)\big)\notag\\
&\in \nabla f^*
\Big(\conv \nabla f\Big( \ffproj{C}\big(\sqrt{a},\sqrt{a}\,\big)\Big)\Big). \notag
\end{align}
In view of the characterization \eqref{e:zchar} of $\bz$, we deduce
that $\bz = \big(\sqrt{a},\sqrt{a}\,\big)$.
\end{proof}

\begin{remark}~\\[-8mm]
\begin{enumerate}
\item The fact that the extreme points $\{\bc_0,\bc_1\}$ play
a role in Example~\ref{ex:halfeuclid} and Example~\ref{ex:kullback}
is not surprising since in these cases $D(\bx,\cdot)$ is convex for every
$\bx\in U$ (see, e.g., \cite{BB01}) so that
\cite[Corollary~32.3.2]{Rock70} applies.
\item
Note that $\bz$ is the arithmetic mean of $\bc_0$ and $\bc_1$
when $D$ is the halved Euclidean distance squared
(Example~\ref{ex:halfeuclid}),
and that $\bz$ is the geometric mean of $\bc_0$ and $\bc_1$
when $D$ is the Kullback-Leibler divergence
(Example~\ref{ex:kullback}).
This might nurture the conjecture that $\bz$ is the
harmonic mean of $\bc_0$ and $\bc_1$ for the
Itakura-Saito distance ---
depending on the location of $a$,
this is sometimes but not always the case
(see Example~\ref{ex:IS} and Lemma~\ref{l:shawn}).
\end{enumerate}
\end{remark}

\begin{example}[Itakura-Saito distance]
\label{ex:IS}
Suppose that $D$ is as in
Example~\ref{ex:distances}\ref{ex:distances:IS}.
Set
\begin{equation}
g = g(a) =
\frac{a(a+1)}{(a-1)^2}\ln\left(\frac{(a+1)^2}{4a}\right)
\quad
\text{and}
\quad
h = h(a) = \frac{2a}{a+1}.
\end{equation}
Then
\begin{equation}
\label{e:japaconcl}
\bz = \begin{cases}
(h,h), &\text{if $g<h$;}\\
(g,g), &\text{if $g\geq h$;}
\end{cases}
\quad
\text{and}
\quad
\ffproj{C}(\bz) =
\begin{cases}
\{\bc_0,\bc_1\}, &\text{if $g<h$;}\\
\{\bc_0,\bc_{1/2},\bc_{1}\}, &\text{if $g\geq h$.}
\end{cases}
\end{equation}
\end{example}
\begin{proof}
Set
\begin{equation}
\bg = (g,g)
\quad
\text{and}
\quad
\bh = (h,h),
\end{equation}
and note that a straightforward computation yields
\begin{equation}
\label{e:keydelta}
\Delta \cap \nabla f^*\Big(\conv\nabla f\big(\{\bc_0,\bc_1\}\big)\Big) =
\{\bh\}.
\end{equation}
Let $\bx=(x,x)\in U\cap \Delta$ and set
\begin{equation}
d_\bx \colon \RR\to[0,\pinf]\colon\lambda\mapsto
D(\bx,c_\lambda).
\end{equation}
Then $\dom d_\bx = \menge{\lambda\in\RR}{\bc_\lambda\in U} =
\left]-1/(a-1), a/(a-1)\right[ \supset [0,1]$.
For every $\lambda\in\dom d_\bx$, we have
\begin{equation}\label{e:itakuradist}
d_\bx(\lambda)
=\ln\left(\frac{(a-1)\lambda+1}{x}\right)+\frac{x}{(a-1)\lambda+1}+\ln\left(\frac{(1-a)\lambda+a}{x}\right)+\frac{x}{(1-a)\lambda+a}-2.
\end{equation}
and hence
\begin{equation}\label{e:itakura1stpartial}
d_\bx'(\lambda)
=\frac{a-1}{(a-1)\lambda+1}-\frac{x(a-1)}{\big((a-1)\lambda+1\big)^2}+\frac{1-a}{(1-a)\lambda+a}-\frac{x(1-a)}{\big((1-a)\lambda+a\big)^2}.
\end{equation}
We note in passing that an elementary calculation results in
\begin{equation}
\label{e:Gexprs}
d_\bx(0)-d_\bx(\thalb)
= \ln\left(\frac{4a}{(a+1)^2}\right)+\left(\frac{(a-1)^2}{a(a+1)}\right)x.
\end{equation}
Now observe that
$d_\bx'(\thalb) = 0$
and that $d_\bx'(\lambda)$ in \eqref{e:itakura1stpartial} is also
a quotient of two polynomials (in $\lambda$), where the numerator is a
polynomial of degree $3$ or less.
Thus, $d_\bx'$ has at most two further roots different from $\thalb$,
which would have to be centered symmetrically around $\thalb$ because of
the symmetry of $d_\bx$ about $\thalb$.
Furthermore, $d_\bx(\lambda)\to\pinf$ as $\lambda$ approaches
either boundary
point of $\dom d_\bx$. Hence, critical points of $d_\bx$ that
are different from
$\thalb$ cannot be local maximizers. Therefore,
$\ffproj{C}(\bx) \subseteq \{\bc_0,\bc_{1/2},\bc_{1}\}$.
The symmetry of $D$ and $C$ yields that exactly one of the following
holds:
$\ffproj{C}(\bx)=\{\bc_{1/2}\}$,
$\ffproj{C}(\bx)=\{\bc_{0},\bc_1\}$,
or
$\ffproj{C}(\bx)=\{\bc_{0},\bc_1,\bc_1\}$.
Combining this with \eqref{e:Gexprs}, we obtain the equivalences
\begin{equation}
\label{e:someequis}
\ffproj{C}(\bx) =\{\bc_{0},\bc_1\}
\quad\Leftrightarrow\quad
d_\bx(0)-d_\bx(\thalb)> 0
\quad\Leftrightarrow\quad
x > g.
\end{equation}
Let us now turn to the
{$\fD{}\,$-Chebyshev center} $\bz$ of $C$.
Since $\bz\in\Delta$ (Proposition~\ref{p:ondiag}) and
$\ffproj{C}(\bz)$ must contain at least 2 points
(Corollary~\ref{c:chebmultval}),
we write $\bz=(z,z)$ and we deduce that
either
$\ffproj{C}(\bz) = \{\bc_0,\bc_1\}$
or
$\ffproj{C}(\bz) = \{\bc_0,\bc_{1/2},\bc_1\}$.
In turn, this means that exactly one of the following two cases holds.
\begin{equation}
\label{e:c1}
\tag{Case~1}
\ffproj{C}(\bz) = \{\bc_0,\bc_1\},
\end{equation}
or
\begin{equation}
\label{e:c2}
\tag{Case~2}
\ffproj{C}(\bz) = \{\bc_0,\bc_{1/2},\bc_1\}.
\end{equation}
If \eqref{e:c1} holds, then
\eqref{e:zchar}, Proposition~\ref{p:ondiag}, \eqref{e:keydelta},
and \eqref{e:someequis} yield
that $\bz = \bh$ and that $z>g$.
Thus,
\begin{equation}
\label{e:c1imp}
\eqref{e:c1}\;\Rightarrow\; z = h > g.
\end{equation}
Using \eqref{e:Gexprs}, we obtain the implication
\begin{equation}
\label{e:c2imp}
\eqref{e:c2}\;\Rightarrow\; z = g.
\end{equation}

We now assume momentarily that $g<h$.
Then, by \eqref{e:someequis}, $\ffproj{C}(\bh)=\{\bc_0,\bc_1\}$ and
hence $\bh\in\nabla f^*(\conv\nabla \ffproj{C}(\bh))$ by
\eqref{e:keydelta}. In view of the characterization \eqref{e:zchar} of
$\bz$, we obtain $\bz=\bh$ and hence $z=h$. We thus have verified the first
case of \eqref{e:japaconcl}.

Finally, we assume that $g\geq h$.
In view of \eqref{e:c1imp}, \eqref{e:c1} cannot hold.
Thus, \eqref{e:c2} must hold and \eqref{e:c2imp} yields that $z=g$, i.e.,
that $\bz=\bg$.
\end{proof}

The formula for $\bz$ given in Example~\ref{ex:IS} immediately
raises the question on how $g$ and $h$ relate to each other, viewed
as functions of $a$.
In the following result, we provide
an alternative description of the inequality
$g<h$.

\begin{lemma}
\label{l:shawn}
Let the functions $g$ and $h$ be defined on the interval
$I = \left]1,\pinf\right[$ by
\begin{equation}
g(x) =
\frac{x(x+1)}{(x-1)^2}\ln\left(\frac{(x+1)^2}{4x}\right)
\quad
\text{and}
\quad
h(x) = \frac{2x}{x+1}.
\end{equation}
Then there exists a
real number $\widetilde{a}\in I$ such that
\begin{equation}\label{howtocompare}
\big(\forall x\in I\big)\quad
\begin{cases}
g(x)<h(x), &\text{if $x<\widetilde{a}$;}\\
g(x)=h(x), &\text{if $x=\widetilde{a}$;}\\
g(x)>h(x), &\text{if $x>\widetilde{a}$.}
\end{cases}
\end{equation}
In fact, $\widetilde{a}\approx 17.63$.
\end{lemma}

\begin{proof}
Observe that
\begin{align}
\label{e:keycharineqs}
h(x)>g(x)
&\Leftrightarrow
\frac{2(x-1)^2}{(x+1)^2}>\ln\bigg(\frac{(x+1)^2}{4x}\bigg)=
2\ln(x+1)-\ln(4x) \\
&\Leftrightarrow
k(x) := \frac{2(x-1)^2}{(x+1)^2} -
2\ln(x+1)+\ln(4x)>0.\notag
\end{align}
Since
\begin{align}
k'(x) &= \frac{8(x-1)}{(x+1)^3} - \frac{2}{x+1} + \frac{1}{x}
= \frac{-(x-1)(x^2-6x+1)}{x(x+1)^3} \\
&=
\frac{-(x-1)\big(x-(3-2\sqrt{2})\big)\big(x-(3+2\sqrt{2})\big)}{x(x+1)^3},
\notag
\end{align}
we set $\xi=3+2\sqrt{2}\approx 5.83$, and we
deduce that $k$ is strictly increasing on
$\left]1,\xi\right]$ and that
$k$ is strictly decreasing on
$\left[\xi,\pinf\right[$.
On the other hand, $k(1)=0$ and $\lim_{x\to\pinf} k(x)=\minf$.
Altogether,
there must exist some number $\widetilde{a}>\xi$ such that
$k>0$ on $\left]1,\widetilde{a}\right[$,
$k(\widetilde{a}) = 0$, and
$k<0$ on $\left]\widetilde{a},\pinf\right[$.
In view of \eqref{e:keycharineqs}, we obtain \eqref{howtocompare}.
Finally, the proclaimed approximation
$\widetilde{a}\approx 17.63$ follows
from Maple, Mathematica, or by simple bisection.
\end{proof}

\begin{remark}
Consider again Example~\ref{ex:IS} and its notation.
Define numbers  $\mu_0,\mu_{1/2},\mu_1$ according to the following two
alternatives:
\begin{equation}
g < h \quad\Rightarrow\quad
\begin{cases}
\mu_0 = \mu_1 = \thalb; \\[+5 mm]
\mu_{1/2} = 0,
\end{cases}
\end{equation}
or
\begin{equation}
g \geq h \quad\Rightarrow\quad
\begin{cases}\mu_0=\mu_1=\displaystyle\frac{(a-1)^2-2a\ln\left(\displaystyle
\frac{(a+1)^2}{4a}\right)}{(a-1)^2\ln\left(\displaystyle\frac{(a+1)^2}{4a}\right)};\\[+13 mm]
\mu_{1/2}=\displaystyle\frac{-2(a-1)^2+(a+1)^2\ln\left(\displaystyle\frac{(a+1)^2}{4a}\right)}{(a-1)^2\ln\left(\displaystyle\frac{(a+1)^2}{4a}\right)}.
\end{cases}
\end{equation}
One may verify that $\{\mu_0,\mu_{1/2},\mu_1\}\subset [0,1]$,
that $\mu_0 + \mu_{1/2} + \mu_1 = 1$, and that
\begin{equation}\label{e:chebcent_muvals}
\bz = \nabla f^*\Big( \mu_0\nabla f(\bc_0) + \mu_{1/2}\nabla f(\bc_{1/2}) +
\mu_{1}\nabla f(\bc_1)\Big).
\end{equation}
Note that the existence of such
convex coefficients is guaranteed by \eqref{e:zchar}.
\end{remark}

\begin{figure}
\centering
\begin{tabular}{cc}
  \begin{tabular}{ccc}
  Euclidean & Kullback-Leibler & Itakura-Saito \\
\includegraphics[width=4cm]{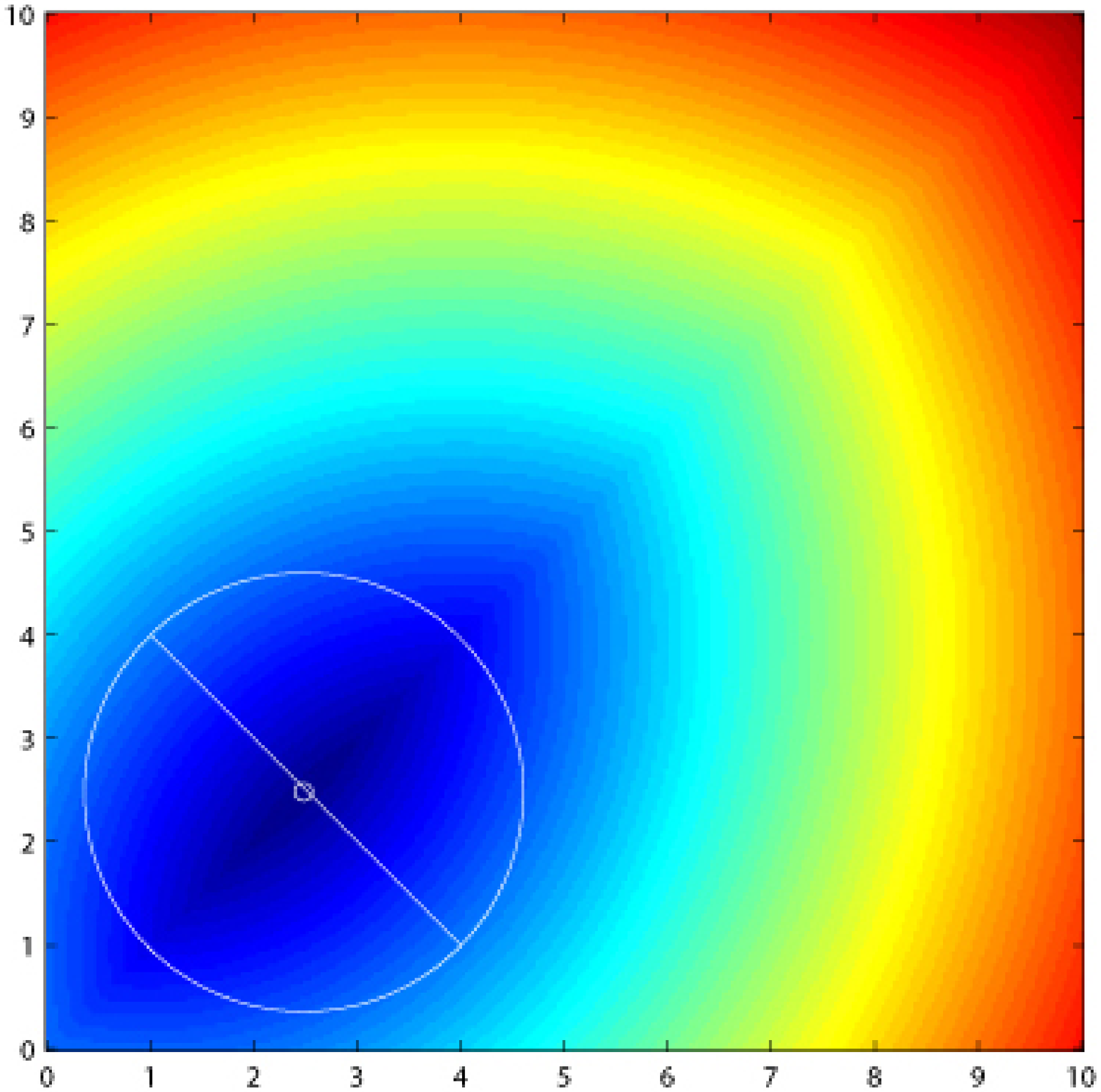} &   \includegraphics[width=4cm]{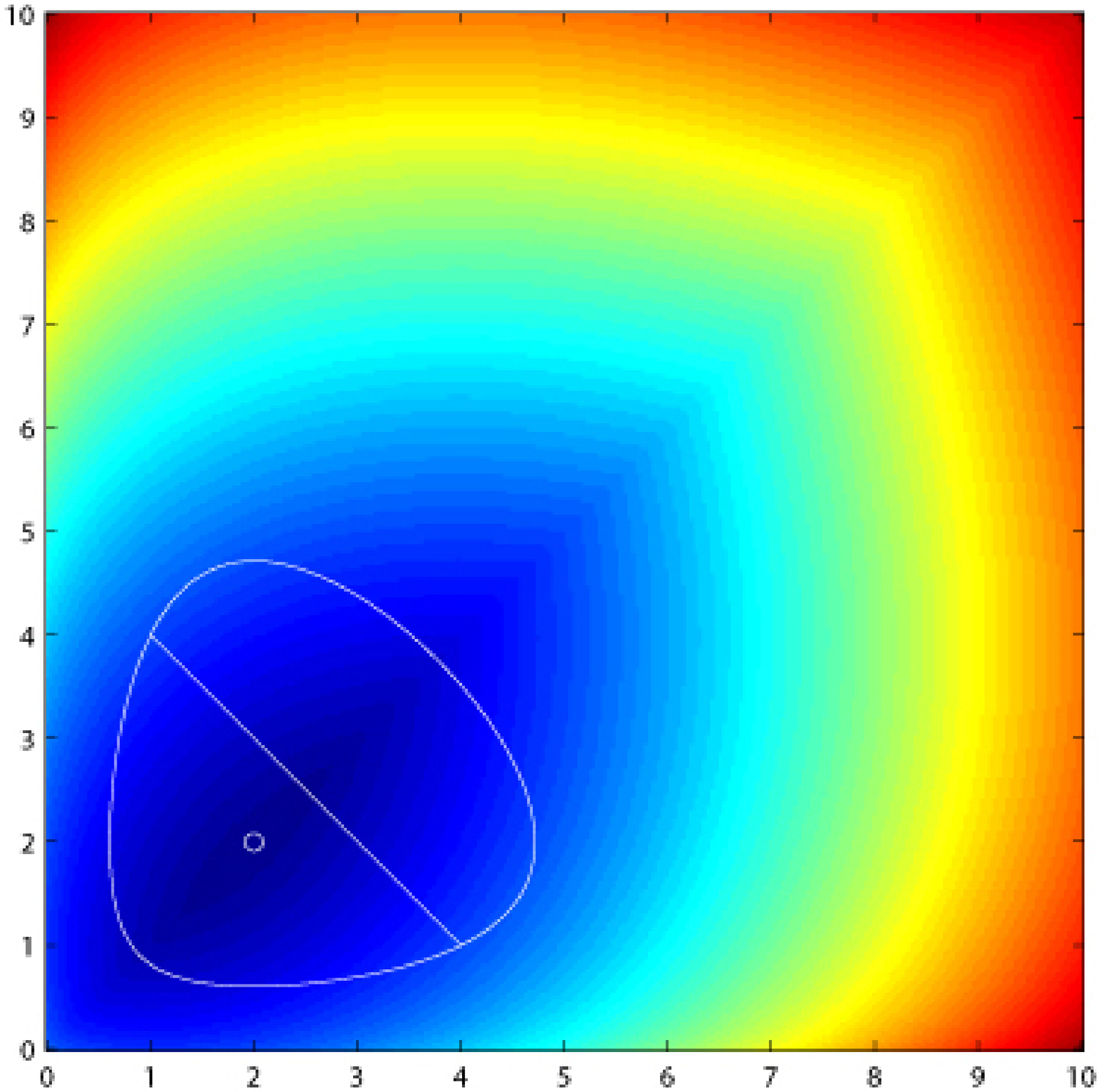} &   \includegraphics[width=4cm]{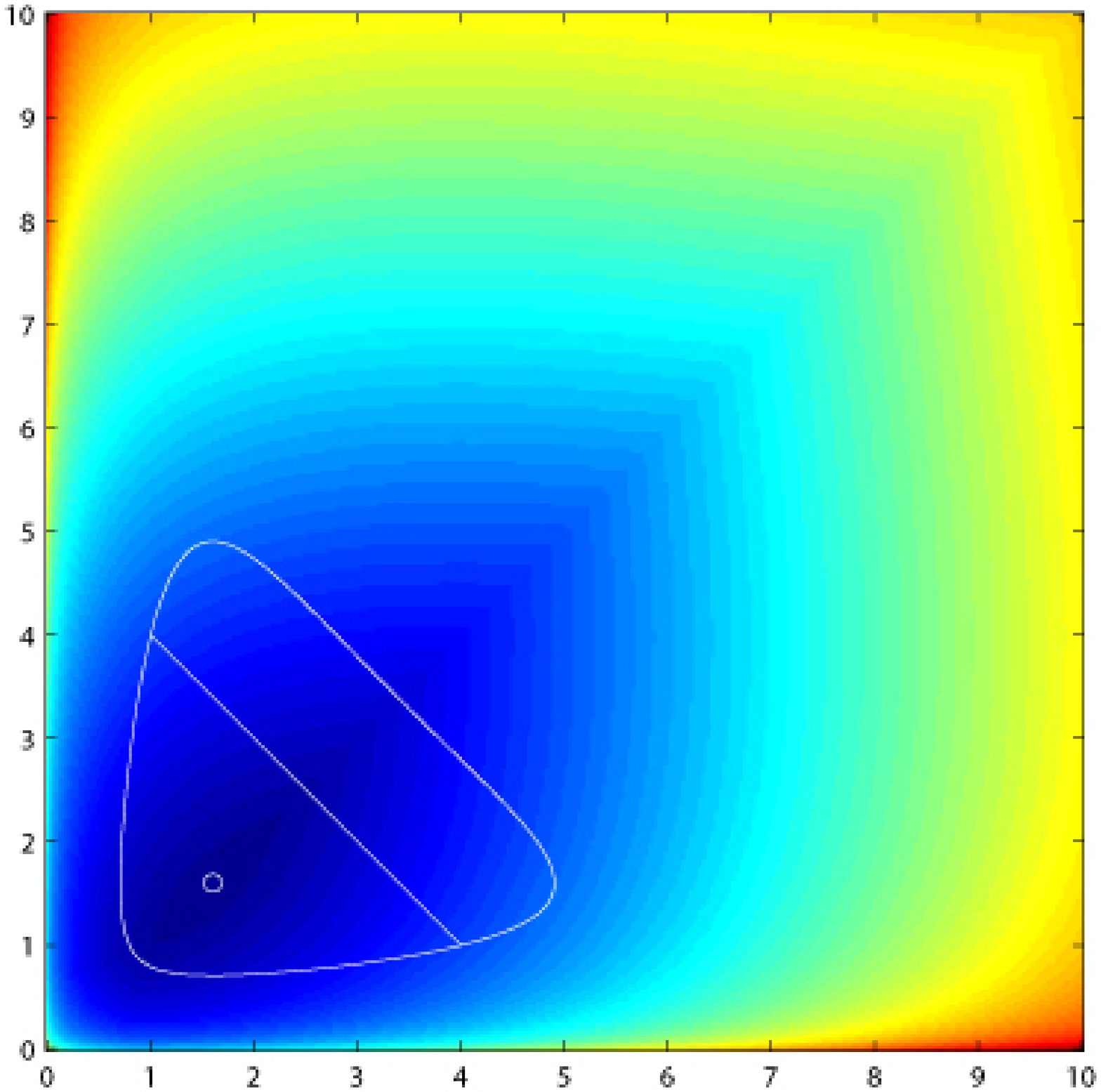}\\
  ~ & $a=4$ & \\
\includegraphics[width=4cm]{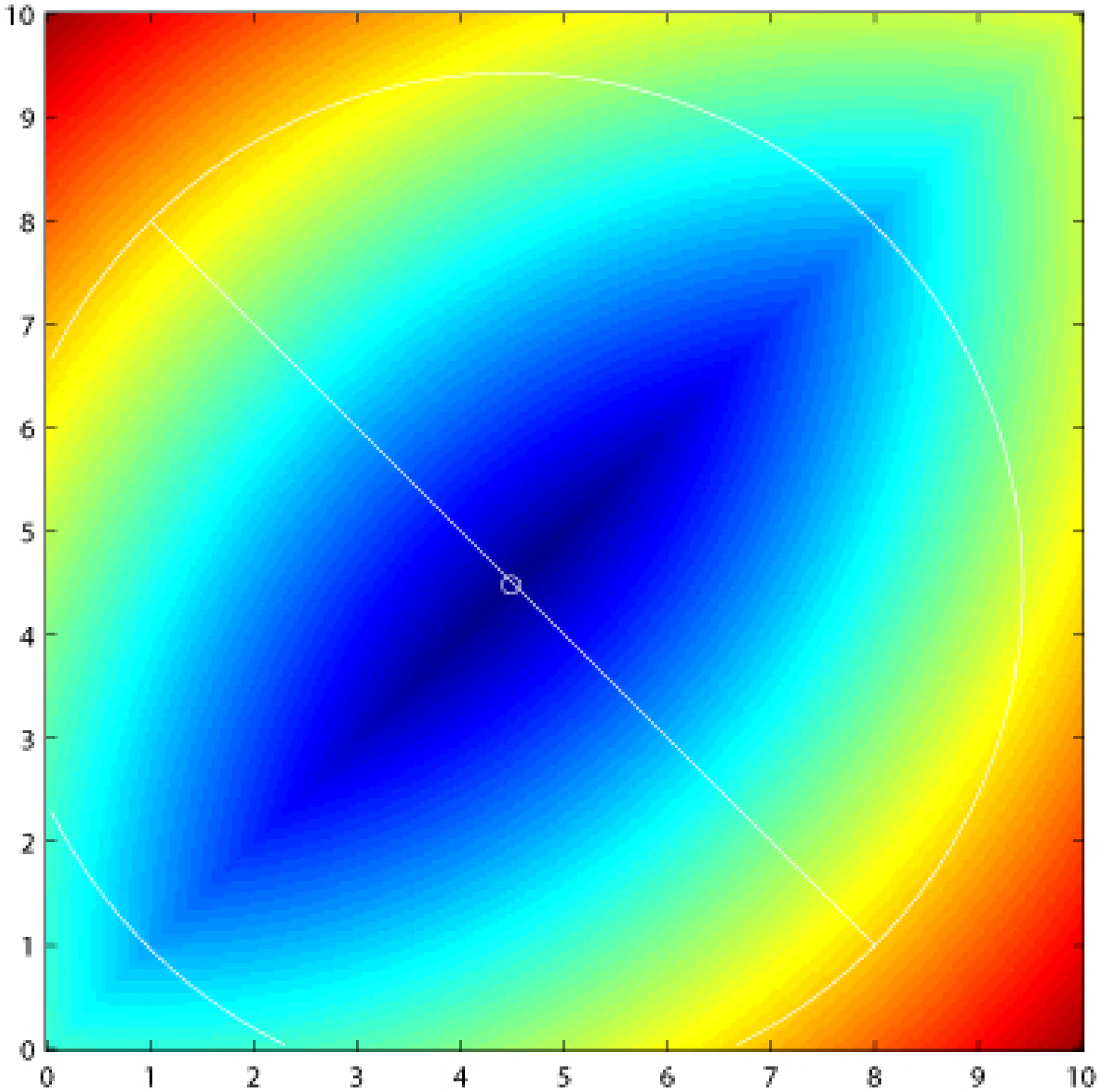} &     \includegraphics[width=4cm]{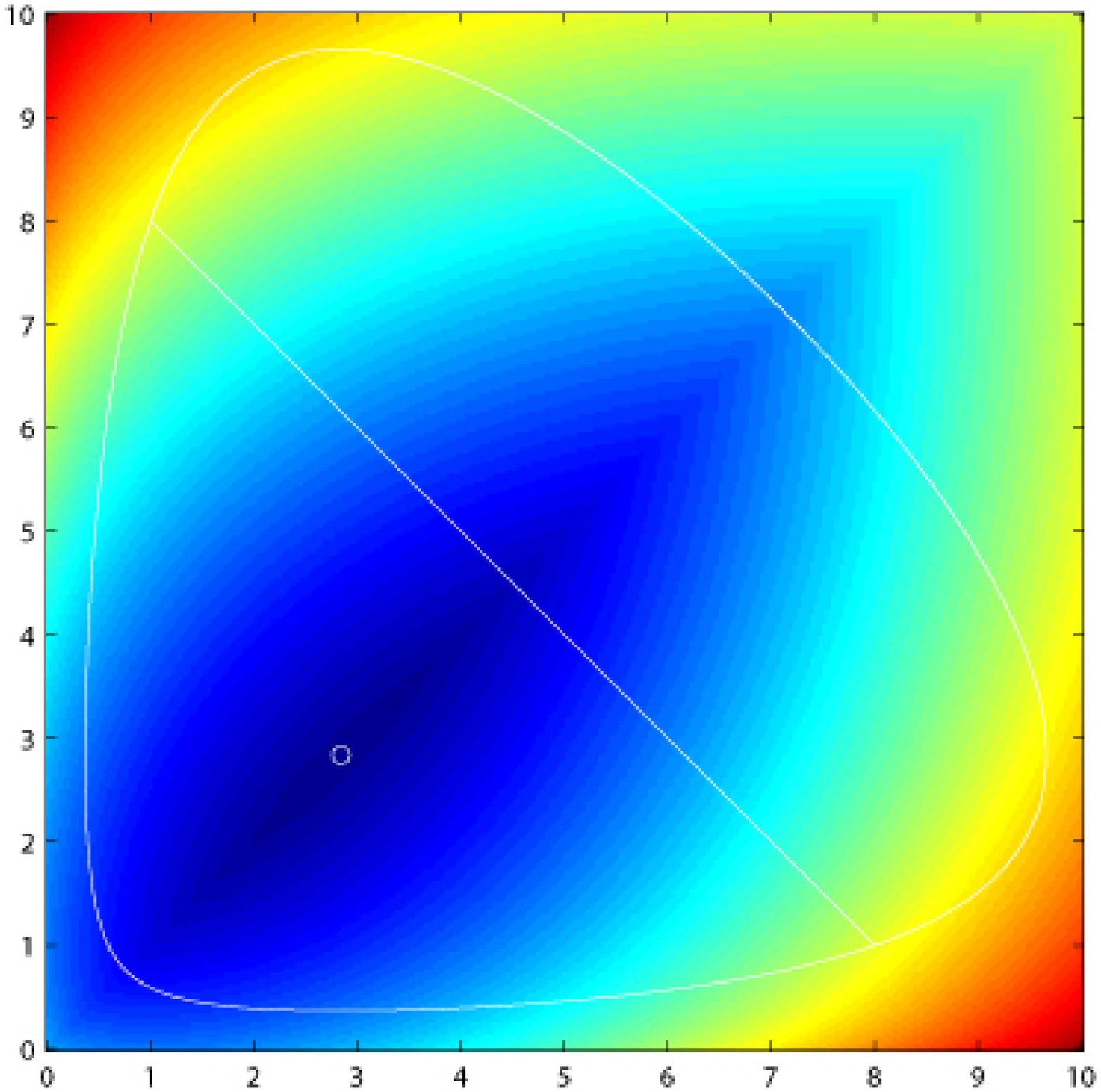} &   \includegraphics[width=4cm]{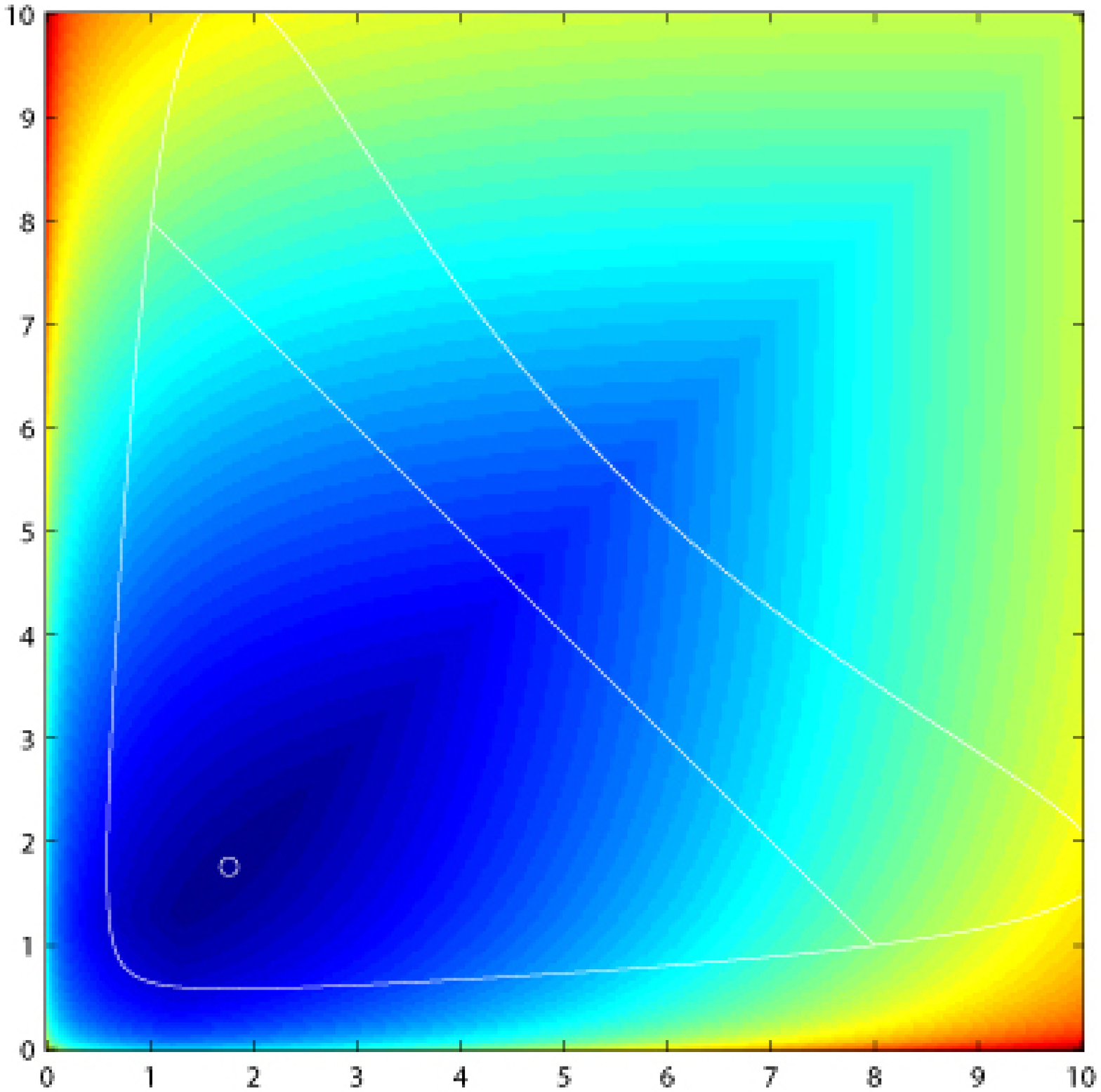}\\
  ~ & $a=8$ & \\
\includegraphics[width=4cm]{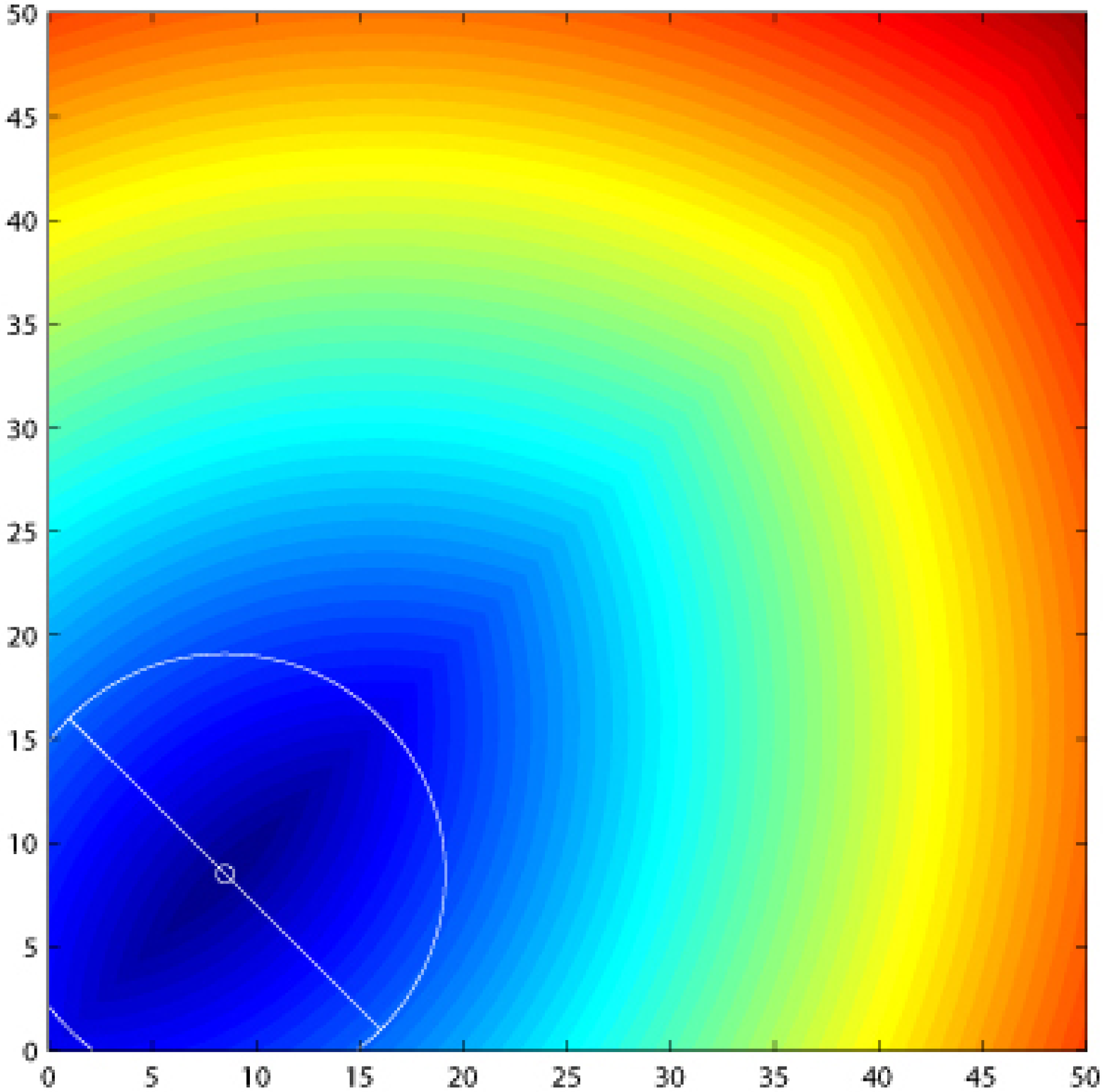} &     \includegraphics[width=4cm]{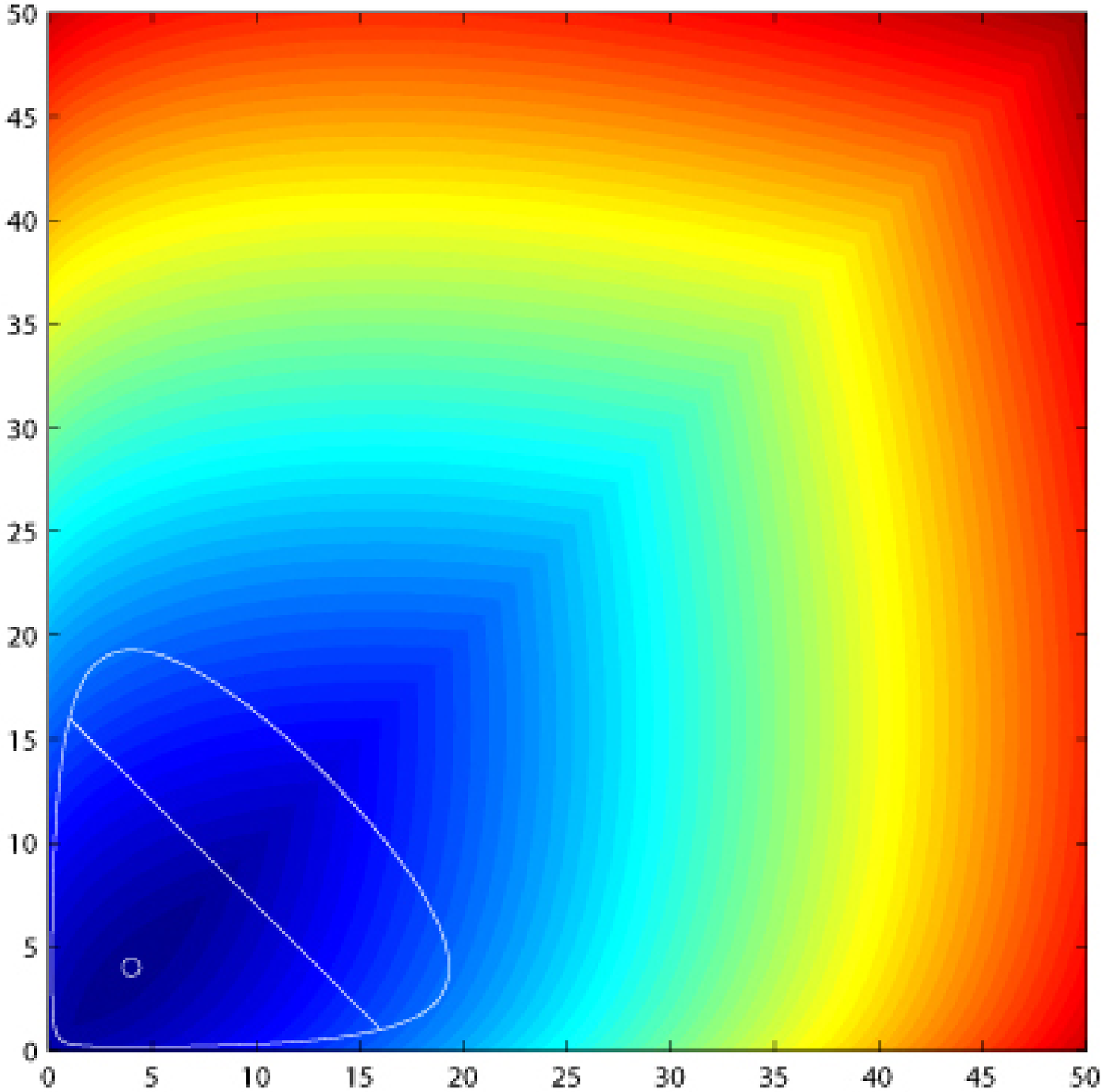} &   \includegraphics[width=4cm]{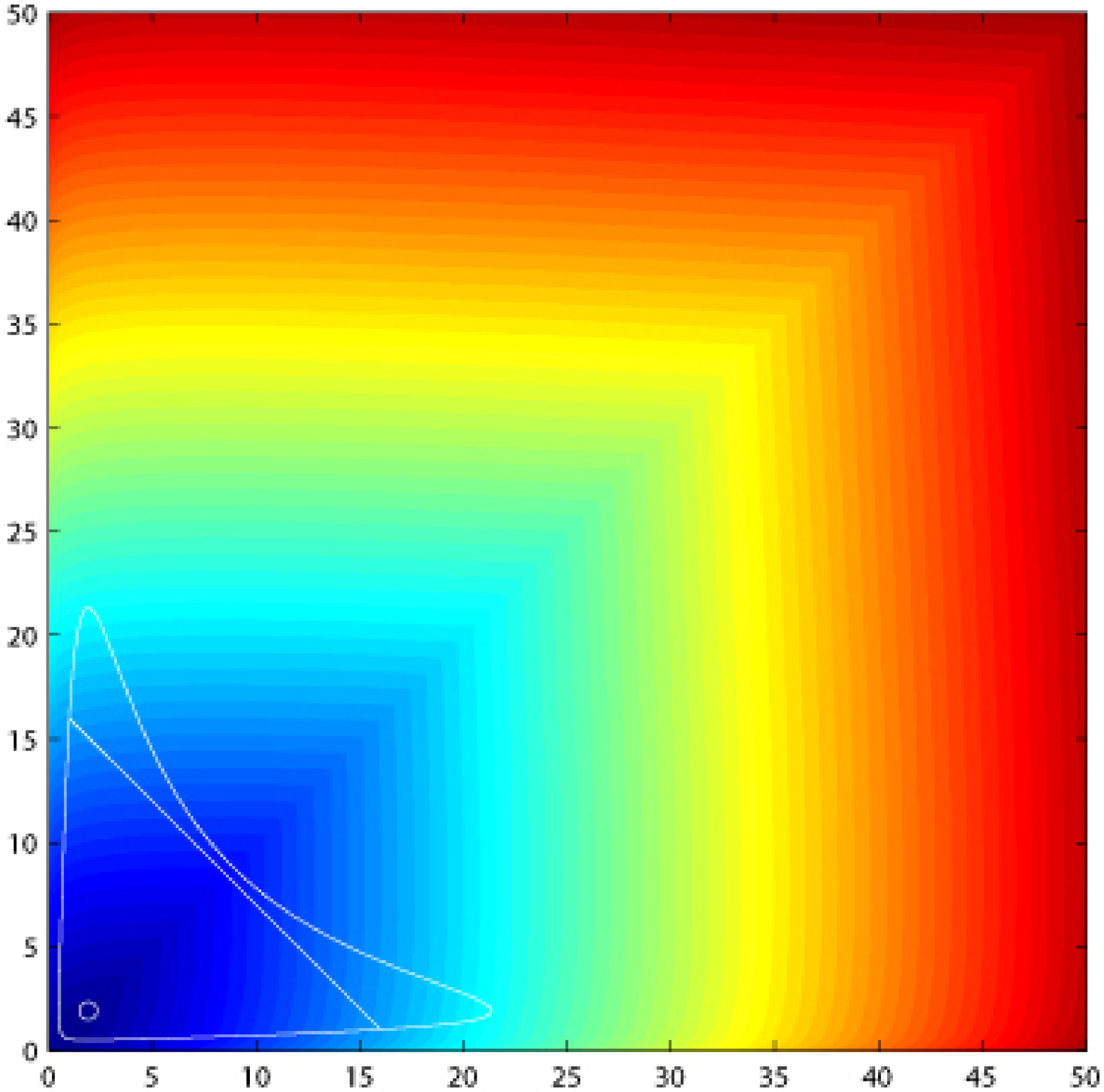}\\
  ~ & $a=16$ & \\
\includegraphics[width=4cm]{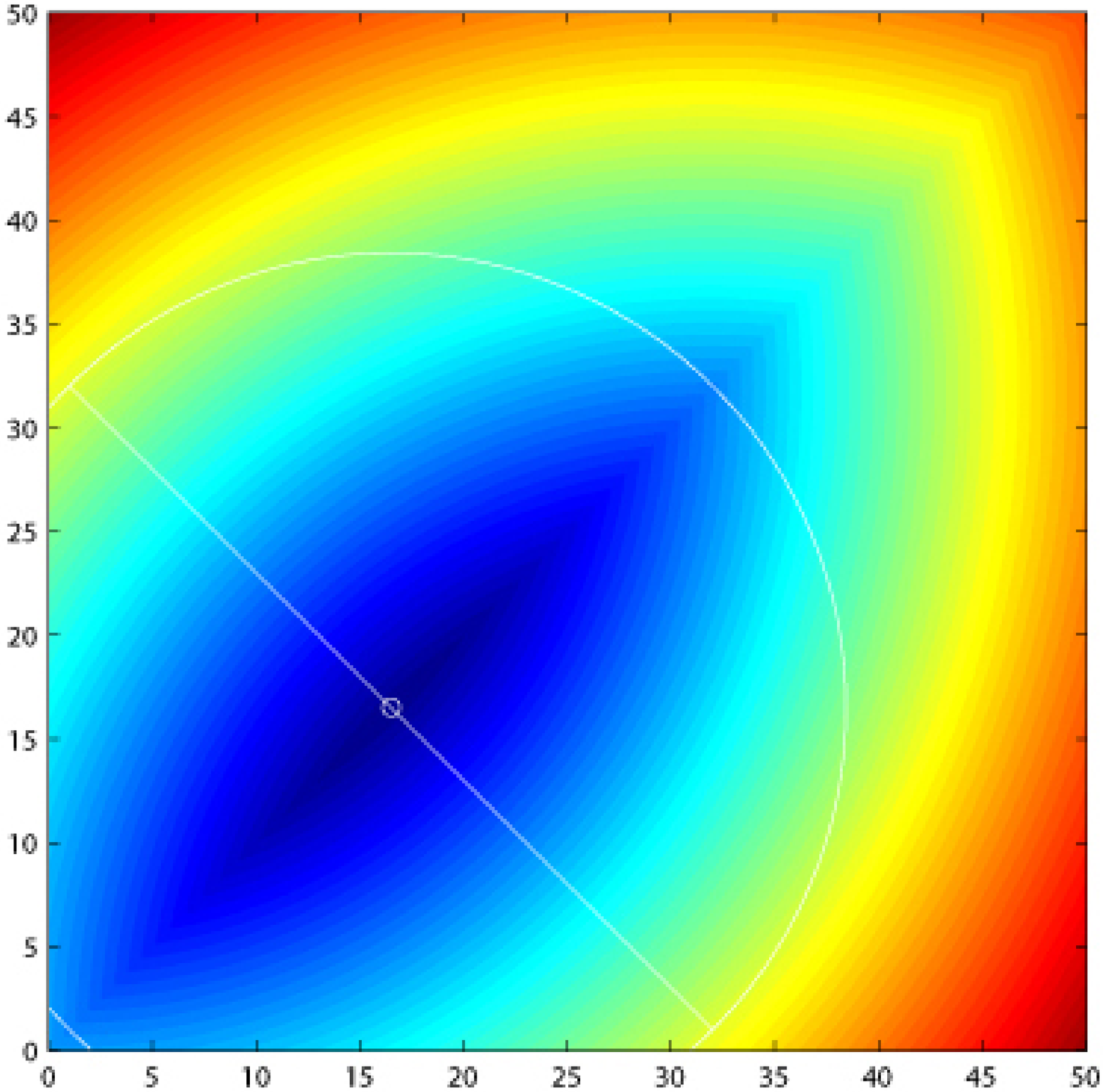} &   \includegraphics[width=4cm]{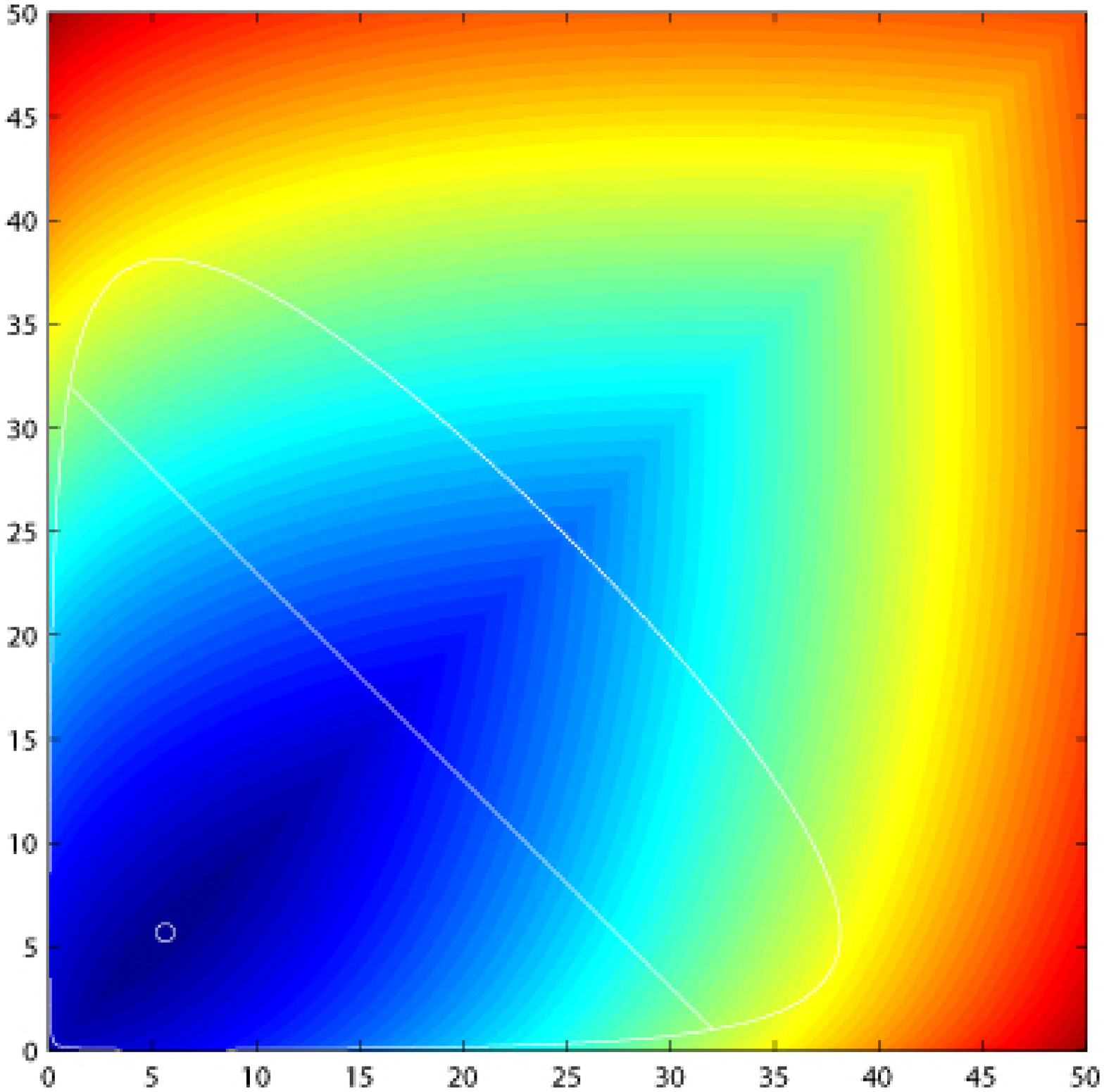} &   \includegraphics[width=4cm]{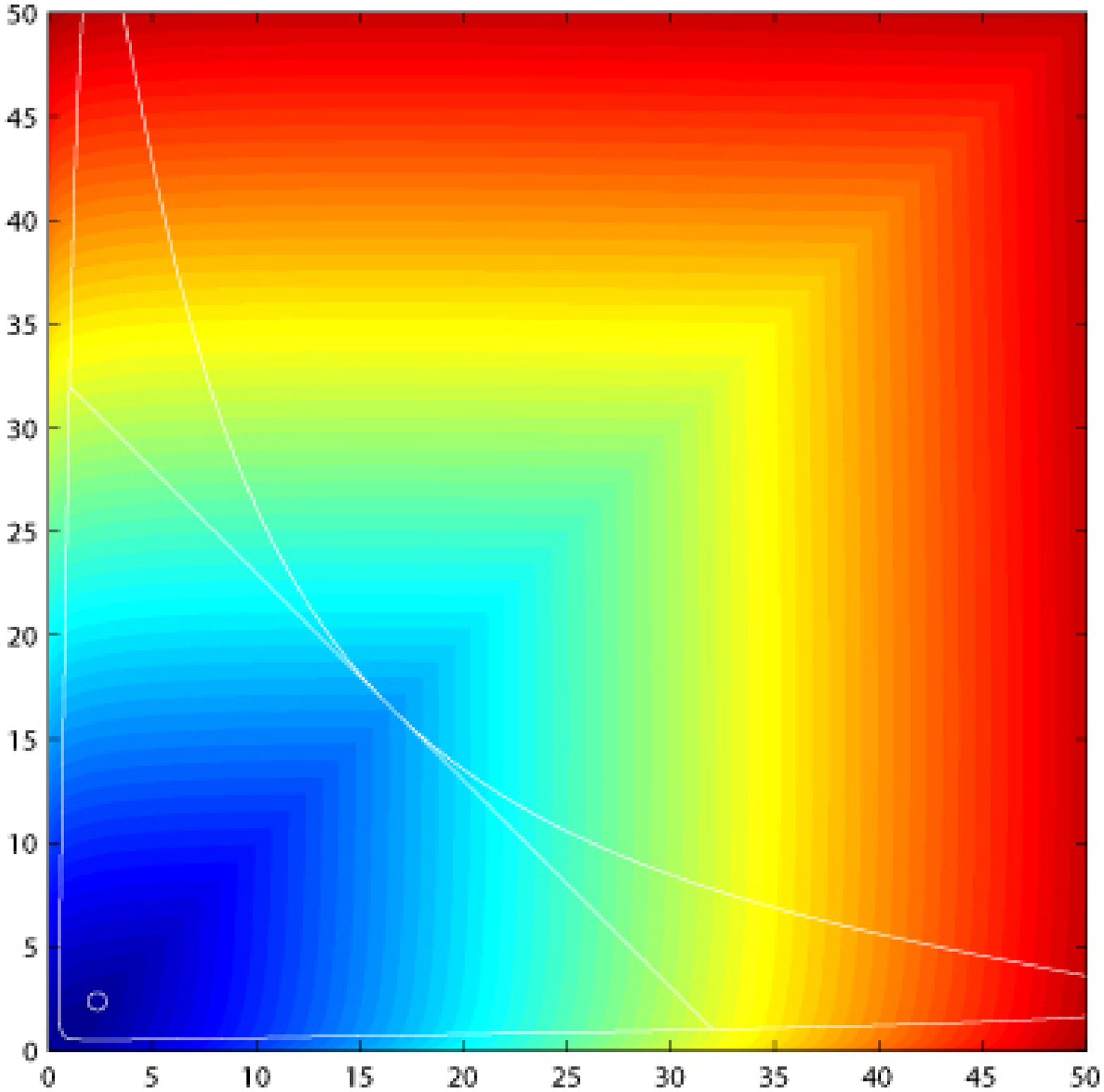}\\
  ~ & $a=32$ & \\
  \end{tabular} &  \begin{tabular}{c}
  \vspace{0.5cm} \\
  \includegraphics[width=1cm]{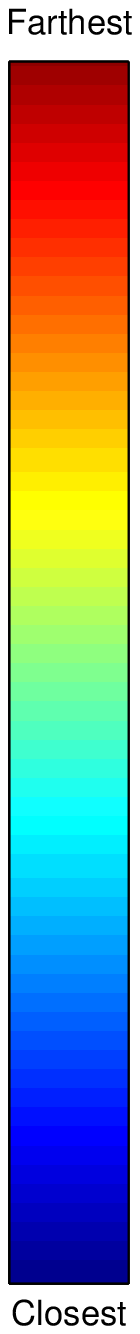}
  \end{tabular}
  \end{tabular}
  \caption{The set $C$, the Chebyshev center $\bz$ of $C$, and the sphere
  of radius $\ffDbz$ centered at $\bz$, for $C$ the line segment connecting
  $(1,a)$ and $(a,1)$ for $a=4$ and $a=8$ over the region
  $[0,10]\times[0,10]$; and $a=16$ and $a=32$, over the region
  $[0,50]\times[0,50]$. Each are shown over color-maps for the three distances analyzed in Section \ref{s:chebcent}, with the interpretations of the colors indicated in the color-legend.}\label{fig:distances_colormaps}
\end{figure}

\begin{remark}\label{r:abarval}
Figure \ref{fig:distances_colormaps} shows the set $C$, the Chebyshev
center $\bz$ of $C$, and the corresponding sphere of radius $\ffD{C}(\bz)$
centered at $\bz$, for a variety of values of $a$ (fixed within each row)
and for each of the three distances analyzed (fixed within each column).
Specifically, shown are $a=4$ and $a=8$ over the region
$R=[0,10]\times[0,10]$ (top two rows), and $a=16$ and $a=32$ over the
region $R=[0,50]\times[0,50]$ (bottom two rows). Each are shown over
color-maps indicating $\ffD{C}(\bx)$ for each $\bx\in R$, with the interpretation of the colors indicated in the accompanying color-legend.
Note that the colors indicate distances from each point in the specified
region to the farthest point in $C$, but are only relative comparisons
within each graph; the same color in separate images does not indicate the
same numerical magnitude, neither for a fixed distance $D$ nor for a fixed
value of $a$. In addition, the color-maps for the halved Euclidean distance
squared and the Kullback-Leibler divergence were calculated using $\ffproj
C(\bx)$ in Examples \ref{ex:halfeuclid} and \ref{ex:kullback},
respectively. However, the color-map for the Itakura-Saito distance was
calculated numerically by a discretization of $C$ due to the absence of
a corresponding formula for $\ffproj C(\bx)$ in Example~\ref{ex:IS}.
We make the following observations directly from Figure \ref{fig:distances_colormaps}:
\begin{enumerate}
\item As predicted by our analysis, for the halved Euclidean distance
squared, $\bz$ falls on the point $\bc_{1/2}$ for all values of $a$
(left-column). The color-map corresponds to $
\max\{D(\bx,\bc_0),D(\bx,\bc_1)\}$, with $D(\bx,\bc_0)=D(\bx,\bc_1)$ along $\Delta$ as per (\ref{e:euclid_farpt}).
\item For the Itakura-Saito distance and for small $a$
(see $a=4$ and $a=8$),
the endpoints $\bc_0$ and $\bc_1$ are the farthest points
from the Chebyshev center $(h,h)$.
When $a\geq \widetilde{a}$ (see Lemma~\ref{l:shawn}), then
the farthest points from $(g,g)$ are $\{\bc_0,\bc_{1/2},\bc_{1}\}$,
and $D((g,g),\bc_{1/2}) < \ffD{C}(h,h)$, visually confirming
that $(g,g)$ is now the Chebyshev center
(see Figure \ref{fig:distances_colormaps} for $a=32$).
\end{enumerate}
\end{remark}

\begin{figure}[h]
\centering
  \includegraphics[width=12cm]{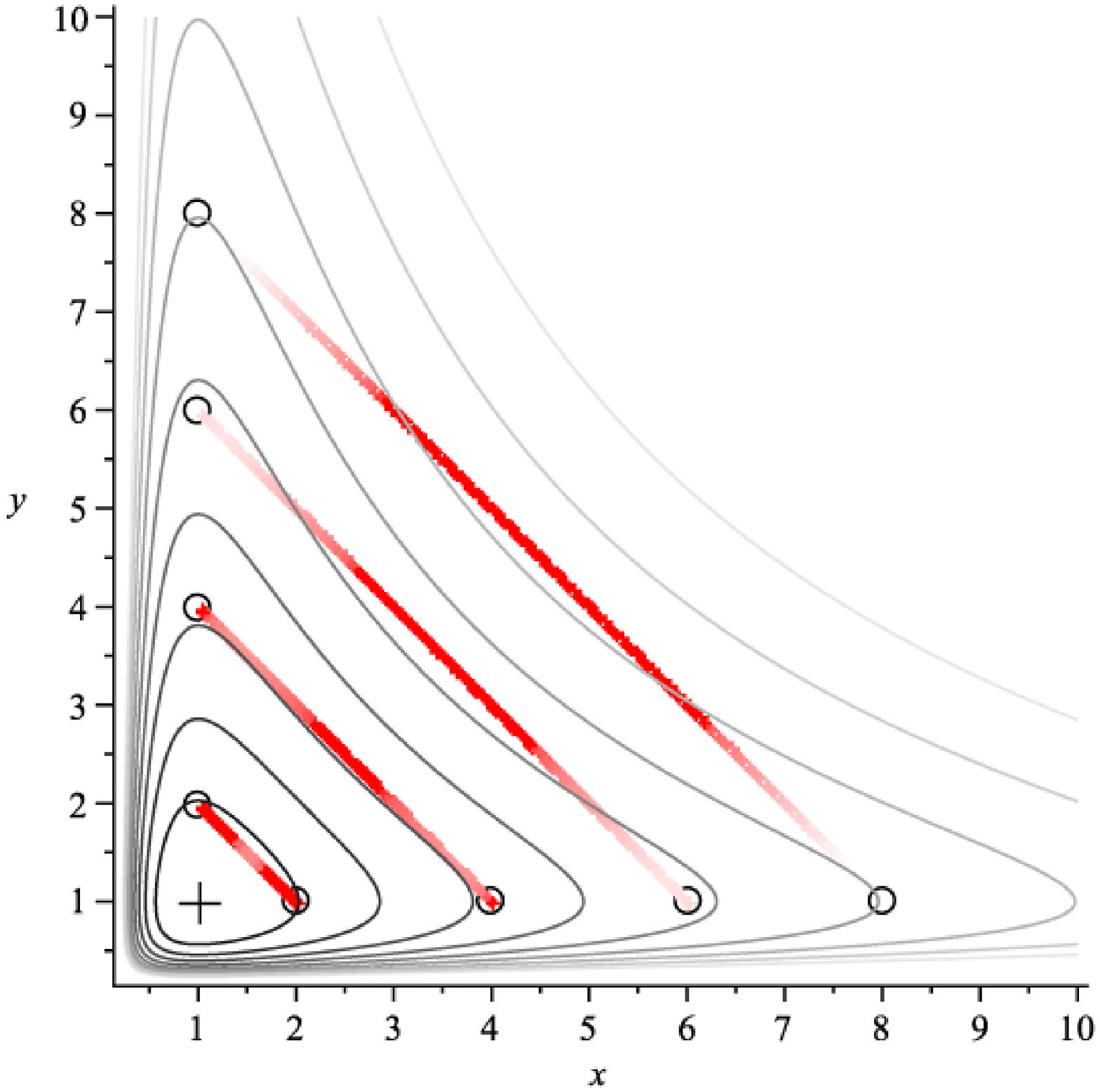}\\
  \caption{Spheres for the Ikaturo-Saito distance centered at $\bx=(1,1)$, for a variety of radii. Also shown are the line-segments $C$ from $(1,a)$ to $(a,1)$ for $a=2,4,6,8$, with plot intensity proportional to the distance from $\bx$.}\label{fig:japanese_linesegment}
\end{figure}

\begin{remark}
Finally, let us fix $\bx=(1,1)$ and assume that $a=6$. For the Itakura-Saito distance, we have that the farthest point $\ffproj C(\bx)$ is $\bc_{1/2}$, which is actually the \emph{nearest} point of $C$ to $\bx$ for both of the other distances. Indeed, Figure \ref{fig:japanese_linesegment} shows the spheres for the Itakura-Saito distance for a variety of radii. The thickness of the line segments is plotted proportional to the distance from $\bx$. (In addition, note that the Itakura-Saito ball is convex for small $a$, a fact not apparent in Figure \ref{fig:distances_colormaps}.)
\end{remark}

\section*{Acknowledgments}
Heinz Bauschke was partially supported by the Natural Sciences and
Engineering Research Council of Canada and by the Canada Research Chair
Program.
Xianfu Wang was partially
supported by the Natural Sciences and Engineering Research Council
of Canada.


\begin{thebibliography}{99}

\bibitem{Edgar1} E.\ Asplund:
``Sets with unique farthest points'',
\emph{Israel Journal of Mathematics},
vol.~5, pp.~201--209, 1967.

\bibitem{Baus97}
H.H.\ Bauschke and J.M.\ Borwein:
``Legendre functions and the method of random Bregman projections'',
\emph{Journal of Convex Analysis}, vol.~4, pp~27--67, 1997.

\bibitem{BB01}
H.H.\ Bauschke and J.M.\ Borwein:
``Joint and separate convexity of the
Bregman distance'' in
\emph{Inherently Parallel Algorithms in
Feasibility and Optimization and
their Applications (Haifa 2000)},
D.\ Butnariu, Y.\ Censor, and S.\ Reich (editors),
pp.~23--36, Elsevier, 2001.

\bibitem{commun01}
H.H.\ Bauschke, J.M.\  Borwein and P.L.\ Combettes:
``Essential smoothness, essential strict convexity, and Legendre
functions in Banach spaces'',
\emph{Communications in Contemporary Mathematics},
vol.~3, pp.~ 615--647, 2001.

\bibitem{BL}
H.H.\ Bauschke and A.S.\ Lewis:
``Dykstra's algorithm with Bregman projections:
a convergence proof'',
\emph{Optimization}, vol.~48, pp.~409--427, 2000.

\bibitem{bwyy1}
H.H.\ Bauschke, X.\ Wang, J.\ Ye and X.\ Yuan:
``Bregman distances and Chebyshev sets'',
\emph{Journal of Approximation Theory},
vol.~159, pp.~3--25, 2009.

\bibitem{bwyy2}
H.H.\ Bauschke, X.\ Wang, J.\ Ye and X.\ Yuan:
``Bregman distances and Klee sets'',
\emph{Journal of Approximation Theory},
vol.~158, pp.~170--183, 2009.


\bibitem{lewis}
J.M.\ Borwein and A.S.\ Lewis,
\emph{Convex Analysis and Nonlinear Optimization},
second edition, Springer-Verlag, 2006.

\bibitem{BorVan01}
J.\ Borwein and J.\ Vanderwerff:
``Convex functions of Legendre type on Banach spaces'',
\emph{Journal of Convex Analysis},
vol.~8, pp.~569--582, 2001.


\bibitem{BorVanBook}
J.\ Borwein and J.\ Vanderwerff,
\emph{Convex Functions: Constructions,
Characterizations \& Counterexamples},
Cambridge University Press, to appear.

\bibitem{Bregman}
L.M.\ Bregman:
``The relaxation method of finding the common point of convex sets and its
application to the solution of problems in convex programming'',
\emph{U.S.S.R.\ Computational Mathematics and Mathematical Physics},
vol.~7, pp.~200--217, 1967.

\bibitem{ButIus}
D.\ Butnariu and A.N.\ Iusem,
\emph{Totally Convex Functions for Fixed Point Computation in Infinite
Dimensional Optimization},
Kluwer, 
2000.

\bibitem{CenZen}
Y.\ Censor and S.A.\ Zenios,
\emph{Parallel Optimization},
Oxford University Press, 1997.

\bibitem{CT}
G.\ Chen and M.\ Teboulle:
``Convergence analysis of a proximal-like minimization
algorithm using Bregman functions'',
\emph{SIAM Journal on Optimization},
vol.~3, pp.~538--543, 1993.


\bibitem{Garkavi}
A.L.\ Garkavi:
``On the \v{C}eby\v{s}ev center and convex hull of a set'' (Russian),
\emph{Uspeshi Matemati\v{c}eskih Nauk}, vol.~19, pp.~139--145, 1964.

\bibitem{urruty2}
J.-B.\ Hiriart-Urruty:
``La conjecture des points les plus \'eloign\'es revisit\'ee'',
\emph{Annales des Sciences Math\'ematiques du Qu\'ebec},
vol.~29, pp.~197--214, 2005.

\bibitem{urruty3}
J.-B.\ Hiriart-Urruty:
``Potpourri of conjectures and open questions in nonlinear analysis
and optimization'',
\emph{SIAM Review},
vol.~49, pp.~255--273, 2007.

\bibitem{Klee60}
V.\ Klee:
``Circumspheres and inner products'',
\emph{Mathematica Scandinavica}, vol.~8, pp.~363--370, 1960.

\bibitem{Klee61}
V.\ Klee:
``Convexity of Chebyshev sets'',
\emph{Mathematische Annalen},
vol.~142, pp.~292--304, 1960/61.

\bibitem{MSV}
T.S.\ Motzkin, E.G.\ Straus, and F.A.\ Valentine:
``The number of farthest points'',
\emph{Pacific Journal of Mathematics},
vol.~3, pp.~221--232, 1953.

\bibitem{NielsenNock}
F.\ Nielsen and R.\ Nock:
``On the smallest enclosing information disk'',
\emph{Information Processing Letters}, vol.~105, pp.~93--97, 2008.

\bibitem{NockNielsen}
R.\ Nock and F.\ Nielsen:
``Fitting the smallest enclosing Bregman ball'',
in \emph{Machine Learning: 16th European Conference on Machine Learning
(Porto 2005)},
J.\ Gama, R.\ Camacho, P.\ Brazdil, A.\ Jorge and L.\ Torgo (editors),
pp.~649--656, Springer Lecture Notes in Computer Science vol.~3720, 2005.

\bibitem{Rock70}
R.T.\ Rockafellar,
\emph{Convex Analysis},
Princeton University Press, Princeton, 1970.

\bibitem{Rock98}
R.T.\ Rockafellar and R. J-B\ Wets,
\emph{Variational Analysis},
Springer-Verlag, 
1998.

\bibitem{Simons}
S.\ Simons,
\emph{From Hahn-Banach to Monotonicity},
Springer-Verlag,
2008.

\bibitem{WesSch}
U.\ Westphal and T.\ Schwartz:
``Farthest points and monotone operators'',
\emph{Bulletin of the Australian Mathematical Society},
vol.~58, pp.~75--92, 1998.

\bibitem{Zali02}{C.\ Z\u{a}linescu},
\emph{Convex Analysis in General Vector Spaces},
World Scientific Publishing, 2002.
\end{thebibliography}
\end{document}